\documentclass[runningheads]{llncs}
\usepackage[T1]{fontenc}
\usepackage{graphicx}
\usepackage{booktabs}
\usepackage[misc]{ifsym}

\usepackage{amssymb}
\usepackage{amsmath}
\usepackage{dsfont}
\usepackage[mathcal]{eucal}
\usepackage{enumitem}
\usepackage{color}
\usepackage{algorithmicx}
\usepackage{algorithm, algpseudocode}

\newcommand{\Dset}{\mathcal{D}}
\newcommand{\Ot}{\tilde{\mathcal{O}}}
\newcommand{\Xs}{{\color{red}X}^\mathrm{s}}
\newcommand{\Ys}{{\color{red}Y}^\mathrm{s}}
\newcommand{\Xv}{{\color{blue}X}^\mathrm{v}}
\newcommand{\Yv}{{\color{blue}Y}^\mathrm{v}}
\newcommand{\ms}{{\color{red}\mu}^\mathrm{s}}
\newcommand{\sigs}{{\color{red}\sigma}^\mathrm{s}}
\newcommand{\sigv}{{\color{blue}\sigma}^\mathrm{v}}
\newcommand{\Ps}{{\color{red}P}^\mathrm{s}}
\newcommand{\Pv}{{\color{blue}P}^\mathrm{v}}
\newcommand{\ks}{{\color{red}k}^\mathrm{s}}
\newcommand{\Ks}{{\color{red}K}^\mathrm{s}}
\newcommand{\Sigs}{{\color{red}\Sigma}^\mathrm{s}}
\newcommand{\Sigv}{{\color{blue}\Sigma}^\mathrm{v}}
\newcommand{\Xe}{X^\mathrm{E}}
\newcommand{\Ye}{Y^\mathrm{E}}
\newcommand{\Pe}{P^\mathrm{E}}

\newtheorem{assumption}{Assumption}

\begin{document}

\title{Time-varying Gaussian Process Bandit Optimization with Experts: no-regret in logarithmically-many side queries}

\titlerunning{Time-varying Gaussian Process Bandit Optimization with Experts}

\author{Eliabelle Mauduit\inst{1} {\Letter} \and
Elo\"ise Berthier\inst{2}  \and
Andrea Simonetto\inst{1}}

\authorrunning{E. Mauduit et al.}

\institute{Unité des Mathématiques Appliquées, ENSTA, Institut Polytechnique de Paris, 91120 Palaiseau, France \email{\{eliabelle.mauduit,andrea.simonetto\}@ensta.fr}
\and
U2IS, ENSTA, Institut Polytechnique de Paris, 91120 Palaiseau, France \email{eloise.berthier@ensta.fr}
}

\tocauthor{Eliabelle~Mauduit,
Elo\"ise~Berthier,
Andrea~Simonetto}

\toctitle{Time-varying Gaussian Process Bandit Optimization with Experts: no-regret in logarithmically-many side queries}

\maketitle              

\begin{abstract}
We study a time-varying Bayesian optimization problem with bandit feedback, where the reward function belongs to a Reproducing Kernel Hilbert Space (RKHS). We approach the problem via an upper-confidence bound Gaussian Process algorithm, which has been proven to yield no-regret in the stationary case. 

The time-varying case is more challenging and no-regret results are out of reach in general in the standard setting. As such, we instead tackle the question of how many additional observations asked to an expert are required to regain a no-regret property. To do so, we formulate the presence of past observation via an uncertainty injection procedure, and we reframe the problem as a heteroscedastic Gaussian Process regression. In addition, to achieve a no-regret result, we discard long outdated observations and replace them with updated (possibly very noisy) ones obtained by asking queries to an external expert. By leveraging and extending sparse inference to the heteroscedastic case, we are able to secure a no-regret result in a challenging time-varying setting with only logarithmically-many side queries per time step. Our method demonstrates that minimal additional information suffices to counteract temporal drift, ensuring efficient optimization despite time variation.

\keywords{Gaussian Processes  \and Upper confidence bounds \and Bandit feedback \and Sparse inference \and Time-varying optimization.}
\end{abstract}

\section{Introduction}

We consider the problem of sequentially optimizing a reward function $f: \Dset \times \mathbb{R}_+ \to \mathbb{R}$ where $\Dset \subset \mathbb{R}^d$ is a compact convex set. In this configuration, the objective depends both on time and on a continuous decision space $\Dset$. At each discrete time step $t$, we obtain a noisy observation of the reward $y_t = f(x_t, t) + \epsilon_t$, where $\epsilon_t \sim \mathcal{N}(0, \sigma^2)$. Our objective is to maximize the sum of rewards
\begin{equation}\label{Eq:MaxReward}
    \underset{(x_t)_t \in \Dset^T}{\text{max}}\hspace{5pt}\sum^T_{t=1} \Big[f(x_t, t) =: f_t(x)\Big].
\end{equation}

At least in the static case, when $f$ does not change in time, this type of problem has often been formulated via Bayesian optimization with bandit feedback~\cite{Srinivas_2012}, whereby an agent must take a sequence of actions while observing the corresponding sequence of rewards. Each action consists of picking a decision $x$ to get an estimate of the reward at the corresponding point. The agent does not modify its environment through its actions and can thus \textit{exploit} previous measurements to predict actions that offer the highest rewards, but should also \textit{explore} new decisions where the value of the reward function is possibly high. For dynamic rewards, the setting is more challenging, as we will see.  

In the time-varying case, the performance metric we are interested in is the dynamic cumulative regret, defined as
\begin{equation}\label{Eq:DynamicRegret}
    R_T = \sum_{t=1}^T \left(\underset{x \in \Dset}{\text{max}} \hspace{5pt} f_t(x) - f_t(x_t)\right)\,,
\end{equation}
representing the cumulative loss in reward picking decision $x_t$ at time $t$ with respect to the best decision \emph{at the same time step}. Algorithms that achieve an asymptotically vanishing average dynamic cumulative regret, as $\lim_{T\to\infty} R_T/T = 0$, are said to enjoy no-regret.     

To derive our main theoretical results, we will work under two reasonable blanket assumptions. First, to model smoothness properties of the functions $f_t$, we assume that they all belong to a Reproducing Kernel Hilbert Space (RKHS) and have bounded RKHS norm. The RKHS associated with kernel $k$ $(\mathcal{H}_k(\Dset),\langle .,.\rangle_k)$ is a subspace of $L_2(\Dset)$ \cite{Scholkopf2005KernelMI} and the associated inner product $\langle .,.\rangle_k$ is such that
$$
\forall f \in \mathcal{H}_k(\Dset) \mbox{, } f(x) = \langle f,k(x,.)\rangle_k\,. 
$$ 
The norm $\|f\|_k$ measures the smoothness of $f$ with respect to the kernel function $k$, therefore assuming $\|f_t\|_k$ is bounded translates into regularity assumptions about the objective.
\begin{assumption}\label{as.1} 
For all time steps, functions $x \mapsto f_t(x)$ belong to a Reproducing Kernel Hilbert Space with continuous bounded kernel $k$ such that $\forall x \in \Dset \mbox{, } k(x,x) \leq M_k^2$ and they have bounded RKHS norm,
\begin{equation}\label{Eq:RKHS}
    \forall t\mbox{, } \|f_t\|_k \leq B.
\end{equation}
\end{assumption}

Second, to model time variation, we assume boundedness of the variations as follows.
\begin{assumption}\label{as.2} 
For the sequence of functions $(f_t)_{t=1}^T$, there exists a bounded constant $\Delta$, such that,
\begin{equation}\label{Eq:TimeVariation}
    \forall t\mbox{, } \sup_{x \in \Dset} |f_{t+1}(x) - f_{t}(x)| \leq \Delta.
\end{equation}
We further let $\Delta =1$, without any loss of generality. 
\end{assumption}

Assumption~\ref{as.2} ensures controlled temporal variations, limiting changes between consecutive iterations and provides a sound framework for uncertainty injection.


\subsection{Related work}


Bayesian optimization in the bandit feedback setting has been studied extensively in the static scenario: the landmark work of Srinivas and coauthors~\cite{Srinivas_2012} proposes an upper-confidence bound algorithm based on a Gaussian Process model of the unknown function obtaining no-regret in several settings. In particular, first, they use the noisy observations of $f$ to derive a possibly miss-specified estimation of its mean $\mu_t(x)$ and covariance $\sigma^2_t(x)$, via a Gaussian Process:
\begin{eqnarray}
\mu_t(x) &=& k_{t}(x)^\top(K_{t}+\Sigma_t)^{-1} Y_t \\
\sigma^2_t(x) &=& k(x,x) - k_{t}(x)^\top(K_{t}+\Sigma_t)^{-1} k_t(x),
\end{eqnarray}
where $\Sigma_t:= \sigma^2 I_t$, $\sigma^2$ being the noise variance of each observation $y_i$, $(K_t)_{i,j}= k(x_i, x_j)$, $k_t(x)= [k(x_1, x), \ldots, k(x_t, x)]$, $k$ being the kernel or covariance function, and $Y_t = [y_1, \ldots, y_t]^\top$. 
Then, since the reward $f$ is unknown, they propose choosing the next decision based on the upper-confidence bound proxy, as,
\begin{equation}
x_{t+1} = \arg\max_{x \in \Dset}\, \mu_t(x) + \beta_{t+1} \sigma_t(x),
\end{equation}
where $(\beta_t)_{t\geq 1}$ is a sequence of positive parameters chosen to ensure a trade-off between exploration and exploitation and it is decisive in proving the convergence of the algorithm. Their algorithm, labeled GP-UCB, obtains no-regret in high probability when $f$ is sampled from a GP, i.e., $f \sim \mathrm{GP}(0, k(x,x'))$ but also for arbitrary $f$ with bounded RKHS norm. As a means of comparison for the square exponential kernel and $f$ having a bounded RHKS norm, they obtain a $R_T = \Ot(\sqrt{T})$ result. Here, the notation $\Ot(\cdot)$ hides poly-logarithmic terms.  

The cited work focused on a noise model whose distribution is identical across observations, also known as homoscedastic setting. Makarova and coauthors in~\cite{DBLP:journals/corr/abs-2111-03637} remove this assumption, define $\Sigma_t:= \mathrm{diag}(\sigma_1^2, \ldots, \sigma_t^2)$ for noise model $\epsilon_t \sim \mathcal{N}(0, \sigma_t^2)$, and deliver a regret bound that matches $R_T = \Ot(\sqrt{T})$ up to a multiplicative $\bar{\sigma} := \max\{\sigma_i\}$ factor, for the heteroscedastic setting.

The time-varying case has also received attention. The work of~\cite{bogunovic2016} extends the GP-UCB algorithm by considering a time-varying reward. They model the time variations by considering a spatio-temporal kernel with a forgetting factor $\varepsilon$, as
\begin{equation}\label{Eq:KernelTV-GP-UCB}
    \forall \hspace{5pt} t_i, t_j \leq t, \hspace{5pt} k((x_{t_i}, t_i), (x_{t_j}, t_j)) \\ =(1-\varepsilon)^{\vert t_i-t_j \vert/2} k(x_i,x_j),
\end{equation}
where $k(\cdot, \cdot)$ is the static kernel. With this modeling, they propose two algorithms: R-GP-UCB runs GP-UCB on windows of size $w \in \mathbb{N}$ and resets at the start of each window. The second one, TV-GP-UCB, uses the spatiotemporal kernel~\eqref{Eq:KernelTV-GP-UCB}. Under this setting, the authors showed that any GP bandit optimization incurs expected regret of at least $\mathbb{E}[R_T]=\Omega(T\varepsilon)$, meaning the algorithm does not enjoy no-regret for fixed $\varepsilon$. This lower bound is not surprising and it also appears in the multi-armed bandit literature~\cite{besbes2014stochastic}. Furthermore, TV-GP-UCB obtains a $R_T = \Ot(T)$ which implies an increasing average cumulative dynamic regret. 

Building on the literature in dynamic (generalized) linear bandits~\cite{Wu_2018,pmlr-v108-zhao20a,russac2021selfconcordant,werge2023bofucb}, in a series of papers~\cite{DBLP:Zhoujournals/corr/abs-2102-06296,deng2022weighted}, new algorithms are proposed in the time-varying setting: a revised R-GP-UCB algorithm, a new sliding-window algorithm SW-GP-UCB, and a weighted algorithm W-GP-UCB. Under the RKHS setting, they either enjoy cumulative dynamic regrets of $\mathcal{O}(T)$ (matching the lower bound), or $\Ot(T)$ for the latter two (with our variation budget expressed in Assumption~\ref{as.2}). The weighted algorithm is interesting, since it starts from a weighted kernel regression,
\begin{equation}
\hat{f} = \arg\min_{f \in \mathcal{H}_k(\Dset)} \sum_{t=1}^T w_t (y_t - f(x_t))^2 + \lambda_t \|f\|^2_k,
\end{equation}
where $\mathcal{H}_k(\Dset)$ is the RKHS on set $\Dset$ and kernel $k$, $w_t$ is a weight, and $\lambda_t \geq 0$ a parameter; they arrive then at the same iterations of Makarova and coauthors in~\cite{DBLP:journals/corr/abs-2111-03637} for the heteroscedastic setting, but with a growing-in-$T$ noise variance.  

Since no-regret is out of reach in the standard setting, the authors of~\cite{imamura2020timevarying} proposed an algorithm capable of dynamically capturing the changes of the objective function, and thereby acquiring more observations when needed. While this does not guarantee no-regret for a constant sampling time, they show an interesting trade-off between sampling and regret. 

Dealing with a spatio-temporal kernel like~\eqref{Eq:KernelTV-GP-UCB} is theoretically challenging. The works of~\cite{van2012Kernel,Brunzema_2022} propose instead to inject uncertainty into old observations. Their starting point is to consider, at every time $t$, that the variance of old observations increases in time (either exponentially or linearly). This is easier to handle since it is now $\Sigma_t$ that changes, but only on the diagonal. Regret results are not provided, but it is not difficult to see that this approach is equivalent to the weighted kernel regression in the RKHS settings and delivers the same $\Ot(T)$ regret. 

In addition to regret analysis, another active field of research in GP regression involves optimization of algorithms complexity. Regression based on GP models becomes impracticable for large datasets as its time complexity scales as $\mathcal{O}\left(N^3\right)$, where $N$ is the number of observations \cite{snelson2005sparse}. The idea of sparse Gaussian Process regression is to approximate the posterior by performing GP regression on a subset of $M \ll N$ inputs. In this way, the complexity becomes $\mathcal{O}\left(NM^2\right)$. The difficulty lies in the selection of the set of sparse inputs (also called pseudo or inducing inputs) and several techniques exist. For example, in \cite{pmlr-v5-titsias09a}, Titsias considers sparse inputs as variational parameters selected to minimize the Kullback-Leibler (KL) divergence between the exact and approximate posteriors. Leveraging this work, Burt et al. show in \cite{JMLR:v21:19-1015} that $M = \mathcal{O}(d\log^d(N))$ sparse inputs suffice to accurately approximate the posterior in terms of KL divergence. They make use of an approximation of a $M$-Determinantal Point Process ($M$-DPP) \cite{Kulesza_2012} to build the set of sparse inputs. $M$-DPPs define a probability distribution over input subsets of size $M$ that favors the selection of dispersed and less correlated points.

Finally, the algorithms developed in the literature show affinity with online learning in the dynamic setting, e.g.,~\cite{pmlr-v38-jadbabaie15}. 

\subsection{Contributions}

In this paper, we extend the literature in several ways. 

\begin{itemize}[label=$\bullet$]
\item First, motivated by the fact that handling time-variations with spatio-temporal kernels is technically challenging, we embrace the uncertainty injection framework and we formulate the time-varying problem as a sequence of static regression problems, with growing-in-time uncertainty. This renders the GP problem a heteroscedastic one.

\item Then, since a no-regret result is out of reach in this setting, we ask \emph{how many additional queries one should pose to an expert} in order to regain the no-regret result that we enjoy in static settings. The answer to the queries are noisy evaluations (or predictions) of the function at a given time. To limit the number of queries, we leverage sparse inference and we estimate the error of updating past observations with the least number of observations as possible.  We call our new GP-UCB algorithm SparQ-GP-UCB for sparse queries. The algorithm performs GP-UCB updates at every time step $t$ by discarding past measurements taken at times $\tau$ farther away than $\mathcal{O}(\log(t))$ steps and asks new observations to an expert. 
 
\item We prove that SparQ-GP-UCB achieves a $R_T=\Ot(\sqrt{T})$ in $\Ot(1)$ additional queries per time step, and it exhibits a $\Ot(T^2)$ computational complexity. This makes SparQ-GP-UCB the first true no-regret time-varying Gaussian Process algorithm, at the expense of logarithmically-many side queries at each step. 
   
\end{itemize}

\section{Problem setting}

\subsection{Uncertainty injection}\label{sec:UI}

We recall our setting. We consider the problem of sequentially optimizing a reward function $f: \Dset \times \mathbb{R}_+ \to \mathbb{R}$ where $\Dset \subset \mathbb{R}^d$ is a compact convex set. In this configuration, the objective depends both on time and on a continuous decision space $\Dset$. At each step $t$, we obtain a noisy observation of the reward $y_t = f(x_t, t) + \epsilon_t$, where $\epsilon_t \sim \mathcal{N}(0, \sigma^2)$. We set $f_t(\cdot) := f(\cdot, t)$ for convenience. 

Our approach of the problem is to inject uncertainty into old measurements and to consider each optimization problem depending on functions $f_t$ as a sequence of separate snapshots.  

Under Assumption~\ref{as.2} on the boundedness of function variations, at time $T$, we can consider that past observations are noisy observations of the current function $f_T$ with zero-mean noise and variance that is increased depending on how old the observations are. In particular, we use independent noise random variables $\epsilon_{t,T}$ and model
\begin{equation}
y_t = f_T(x_t) + \epsilon_{t,T}, \qquad \epsilon_{t,T} \sim \mathcal{N}(0, \sigma^2 ([T-t]^2+1)), \quad t\leq T,
\end{equation}
that is the noise standard deviation increases linearly in time. This is similar to the approach of~\cite{Brunzema_2022} involving a Wiener process and it is well-motivated by the fact that the maximum variation of the function between $t$ and $T$ is $T-t$. In fact, if we model $f_{t+1}(x) = f_t(x) + v$, with $v \in \mathcal{U}(-1, 1)$, i.e., the uniform distribution on $[-1, 1]$, then the observation $y_t$ of function $f_t$ can be interpreted as an observation of function $f_T$ with noise $\mathbb{E}[\epsilon_{t,T}] = \mathbb{E}[(T-t)v + \epsilon_t] =0$ and variance $\mathbb{E}[\|\epsilon_{t,T}\|^2] = \sigma^2 (\frac{1}{3\sigma^2}[T-t]^2+1)$. The latter justifies the expression of the noise, up to asymptotically-unimportant constants. 



At time $t$, then, we would like to maximize the reward of $f_t(x)$ by choosing the next action based on past observation $Y_t = [y_1, \ldots, y_t]^\top$, each with its own zero-mean noise and time-dependent variance. We approach this as a heteroscedastic Gaussian Process and perform the update,
\begin{eqnarray}
\mu_t(x) &=& k_{t}(x)^\top(K_{t}+\Sigma_t)^{-1} Y_t \\
\sigma^2_t(x) &=& k(x,x) - k_{t}(x)^\top(K_{t}+\Sigma_t)^{-1} k_t(x),
\end{eqnarray}
where $\Sigma_t:= \mathrm{diag}(\text{Var}(\epsilon_{1,t}), \ldots, \text{Var}(\epsilon_{t,t})=\sigma^2)$, the kernel matrix $(K_t)_{i,j}= k(x_i, x_j)$, and $k_t(x)= [k(x_1, x), \ldots, k(x_t, x)]$.
We choose the next decision as,
\begin{equation}
x_{t+1} = \arg\max_{x \in \Dset}\, \mu_t(x) + \beta_{t+1} \sigma_t(x),
\end{equation}
where $(\beta_t)_{t\geq 1}$ is a sequence of parameters chosen to ensure a trade-off between exploration and exploitation.

As said, a basic version of this update would lead an increasing average regret. To limit the regret, we consider only recent observations and summarize and update the remaining ones. 

\subsection{Sparse inference}

To summarize and update past observations, we leverage and extend recent results from sparse inference provided in~\cite{JMLR:v21:19-1015}. Consider $Y_T$ observations performed at $X_T = [x_1, \ldots, x_T]$ points, as well as the mean and variance function coming from a GP regression on these points. Burt and coauthors in~\cite{JMLR:v21:19-1015} offer an algorithm to select $\Ot(1)$ points in the domain $\Dset$ which would deliver the same mean and variance up to a tunable multiplicative error term. We summarize their main result in the following proposition.

\begin{proposition}[\cite{JMLR:v21:19-1015}]\label{Prop:PosteriorError}
    Consider the problem of estimating an unknown function $f: \Dset \to \mathbb{R}$ via $T$ noisy observations, $y_t = f(x_t) + \epsilon_t$, $\epsilon_t \sim \mathcal{N}(0, \sigma^2)$ acquired at i.i.d. training inputs $X_T$. Let $f$ be a sample path of a Gaussian Process with zero mean and kernel $k$. Consider a squared exponential kernel function $k$ for simplicity. Let $\mu_0(x)$ and $\sigma^2_0(x)$ represent the mean and variance of the Gaussian Process regression performed on the observations. 
    
    Select a tolerance level $\eta\leq 1/5$. Then, there exists an algorithm that selects $\Ot(1) < T$ points in the domain $\Dset$ and their observations $y_t = f(x_t) + \epsilon_t$, $\epsilon_t \sim \mathcal{N}(0, \sigma^2)$, such that if we let $\mu_1(x)$ and $\sigma^2_1(x)$ represent the mean and variance of the Gaussian Process regression performed on the new points and observations, we obtain in high probability,
    \begin{equation}
        \begin{split}
            \vert \mu_1 - \mu_0\vert \leq \sigma_0 & \sqrt{\eta} \leq \frac{\sigma_1 \sqrt{\eta}}{\sqrt{1-\sqrt{3\eta}}},\\
            \vert 1 - \sigma^2_1 / \sigma_0^2 \vert &\leq \sqrt{3\eta}.
        \end{split}
    \end{equation}
\end{proposition}

Proposition~\ref{Prop:PosteriorError} is a condensed version of Proposition~1, Theorem 14, and Corollary 22 in~\cite{JMLR:v21:19-1015}.  

A possible algorithm proposed in the paper to determine the sparse inputs is an approximate determinantal point process (DPP). Such algorithm selects $M < T$ points in order to minimize the difference between the KL divergence of the exact posterior and approximated one. Specifically, for $\epsilon > 0$, one can use an MCMC algorithm, as specified in Algorithm~1 of \cite{JMLR:v21:19-1015} to obtain an $\epsilon$ approximation of a $M$-DPP, with $T$ inputs, with a computational complexity that is upper bounded by $\mathcal{O}\left(T M^3(\log\log T+\log M+\log 1/\epsilon^2)\right)$. 

The most important feature of the DPP algorithm in~\cite{JMLR:v21:19-1015} and Proposition~\ref{Prop:PosteriorError} is that these results do not depend on the observations values $Y_T$, but only on the points $X_T$ where these observations are taken. We will see next how this is key in devising our sparse algorithm and, along the way, how we can extend Proposition~\ref{Prop:PosteriorError} to the heteroscedastic and deterministic setting. 


\section{SparQ-GP-UCB Algorithm}\label{sec:algo}

With all the previous preliminaries in place, we are now ready for the main algorithm: SparQ-GP-UCB. 

The algorithm works in rounds. At each time step $t$, we consider the problem of maximizing the regret $f_t$, with observations $Y_t = [y_1, \ldots, y_t]^\top$ at points $X_t = [x_1, \ldots, x_t]$. The observations are properly injected with uncertainty, so that their variance grows in time as,
\begin{equation}
\epsilon_{i,t} \sim \mathcal{N}(0, \sigma^2 ([t-i]^2+1)), \quad i\leq t.
\end{equation}

The {\bf first} step of the algorithm is to discard observations that have variance greater than $g(t)$, where $g(t) = o\left(t^{1/4}\right)$. We take $g : t \mapsto \sigma^2 \log(t)$ as an illustrative example but any function $g : t \mapsto g(t) = o\left( t^{1/4} \right)$ would work with no change in the proof arguments (and we further discuss it in the proof).

{\bf Second}, we act as if we had access to updated noisy observations for the discarded measurements, with noise being zero-mean and with $\bar{\sigma}^2$ variance. With this pretend observations and the most recent ones with noise less than $\sigma^2\log(t)$, we perform sparse variational inference. We use the approximate DPP algorithm in \cite{JMLR:v21:19-1015} (Algorithm~1) to find the locations $X^{\mathrm{E}} = [x_1^{\mathrm{E}}, \ldots]$, with $|X^{\mathrm{E}}| = \Ot(1) \ll t$ at which to ask an expert for noisy updated observation with zero-mean and variance $\bar{\sigma}^2$. The new expert-delivered observations, together with the most recent ones are guaranteed to be a good approximation of the pretend setting. 

{\bf Third}, we let $\Ys_t, \Xs_t$ being the set of expert-delivered observations together with the most recent ones with noise less than $\sigma^2\log(t)$ and the points at which they are taken. With this, we can compute the mean and variance as,
\begin{eqnarray}
\ms_t(x) &=& \ks_{t}(x)^\top(\Ks_{t}+\Sigs_t)^{-1} \Ys_t \label{eq.sbu1}\\
(\sigs)^2_t(x) &=& k(x,x) - \ks_{t}(x)^\top(\Ks_{t}+\Sigs_t)^{-1} \ks_t(x), \label{eq.sbu2}
\end{eqnarray}
where $\Sigs_t$ is a diagonal matrix containing all the observation variances up to $\max\{\bar{\sigma}^2,\sigma^2 \log(t)\}$, and the kernel elements $\Ks_t$, $\ks_t$, are evaluated on $\Xs_t$. 

And finally, we compute the next decision, via the UCB proxy:
\begin{equation}\label{ucb.s}
x_{t+1} = \arg\max_{x \in \Dset}\, \ms_t(x) + \beta_{t+1} \sigs_t(x).
\end{equation}

The algorithm is summarized in Algorithm~\ref{alg:SparQ-GP-UCB}. We remark the need for performing sparse inference based on Algorithm~1 of \cite{JMLR:v21:19-1015}, whose details are reported in the Appendix. 

\begin{algorithm}[H]
\caption{SparQ-GP-UCB}\label{alg:SparQ-GP-UCB}
\begin{algorithmic}[1]
\Require Domain $\Dset$, kernel $k$
\For{$t=1,2,\hdots$}
	\State Sample $y_t = f_t(x_t) + \epsilon_{t}$
	\State Discard all the observations with a noise $> \sigma^2 \log(t)$
	\State Perform sparse inference on $X_t$ to obtain locations $X^{\mathrm{E}}$ of cardinality $\Ot(1)$
	\State Query an expert to obtain updated observations on $X^{\mathrm{E}}$ for $f_t$
	\State Perform Bayesian updates~\eqref{eq.sbu1}-\eqref{eq.sbu2} to obtain $\ms_t$ and $\sigs$ using $(\Xs_t, \Ys_t)$
    \State Choose the next action $x_{t+1}$ via~\eqref{ucb.s}
\EndFor
\end{algorithmic}
\end{algorithm}

\subsection{Main results}

In this subsection, we report the main results for our algorithm. They are both given for a squared exponential kernel for simplicity, but they can easily be extended to other standard kernels (Matérn for example). 

\begin{theorem}\label{Thm:regret}{(Regret bound for SparQ-GP-UCB)}
    Take any $0 < \delta \leq 1$ and consider a sequence of reward functions $(f_t)_t$ and the observations $y_t = f_t(x_t) + \epsilon_t$, for $\epsilon_t\sim\mathcal{N}(0, \sigma^2)$ i.i.d.. Let Assumptions~\ref{as.1} and \ref{as.2} hold and consider a squared exponential kernel $k$. Let $T$ be a time horizon and $(x_t)_{t=1}^T$ the set of actions chosen by SparQ-GP-UCB (Algorithm \ref{alg:SparQ-GP-UCB}) and set $(\beta_t)_{t=1}^T$ as
$$    
\beta_t =  \sqrt{2\log \left( \frac{2\vert \textcolor{red}\Sigma^\mathrm{s}_{t} + \textcolor{red}{K}^\mathrm{s}_{tt} \vert^{1/2}} {\delta\vert \textcolor{red}\Sigma^\mathrm{s}_{t}\vert^{1/2}} \right)} + \|f_t\|_k .  
$$    
Then, with probability at least $1-\delta$, by asking to an expert $\mathcal{O}\left(\log^d (t)\right)$ queries per time step, SparQ-GP-UCB attains a cumulative dynamic regret of 
    \begin{equation}\label{Eq:RegretBound}
    R_T = \mathcal{O}\left(\sqrt{ Td \log^{d+5}(T)}\sqrt{\log\left( \frac{1}{\delta} \right) +  d \log^{d+3}\left( d\log\left(T \right) \right)}\right) = \Ot(\sqrt{T}).
    \end{equation}
\end{theorem}

The asymptotic bound given in Eq.\eqref{Eq:RegretBound} implies no-regret with probability $1-\delta$ and matches the static case up to poly-logarithmic factors. The main steps of the proof are given in Section~\ref{sec:proofs}. 

We can now discuss briefly the impact of Assumption~\ref{as.2} for the regret bound. A common metric in the literature to account for the time-varying nature of the objective is the variation budget $V_T$ \cite{pmlr-v38-jadbabaie15,deng2022weighted,besbes2014stochastic} defined as,
\begin{equation}\label{Eq:VarBudget}
    \forall T \mbox{, } V_T = \sum_{t=1}^{T-1} \|f_{t+1}-f_t\|_k.
\end{equation}
Let $x \in \Dset$. Then, by the reproducing property of RKHS,
\begin{equation}\label{Eq:BoundSqDiff}
        \vert f_{t+1}(x)-f_t(x)\vert = \vert\langle f_{t+1}- f_t, k(x,.) \rangle \vert
         \leq \|f_{t+1}-f_t\|_k \|k(x,.)\|_k,
    \end{equation}
where we applied Cauchy Schwarz in the RKHS to obtain the inequality. From the reproducing property, $\|k(x,.)\|_k^2=k(x,x)$. As we are working with bounded kernels (Assumption \ref{as.1}) and $\forall x \in \Dset \mbox{, } k(x,x) \leq M_k^2$, we take the infinite norm in the left side of Eq. \eqref{Eq:BoundSqDiff} and sum over $t=1$ to $T-1$ to obtain:
\begin{equation}\label{Eq:IneqSum}
    \sum_{t=1}^{T-1} \|f_{t+1}-f_t\|_{\infty} \leq  M_kV_T\,.
\end{equation}
By Assumption~\ref{as.2}, SparQ-GP-UCB does work even in the case of $\sum_{t=1}^{T-1} \|f_{t+1}-f_t\|_{\infty} = (T-1)$, meaning that our algorithm can achieve no-regret even for a variation budget that grows linearly in time. This improves the result of \cite{deng2022weighted} that requires $V_T = o(T)$ when $V_T$ is known and $V_T = o(T^{1/4})$ otherwise to obtain sublinear regret.

Along with a no-regret result, we also provide a computational complexity estimate as follows.

\begin{theorem}\label{Thm:complexity}
Under the same setting of Theorem~\ref{Thm:regret}, the computational complexity of SparQ-GP-UCB is upper bounded by $\mathcal{O}\left(T^2 \log(T) \log^{3d}\left(\frac{T}{\log(T)}\right)\right) = \Ot(T^2)$.
\end{theorem}

The theorem shows how SparQ-GP-UCB is actually less computationally expensive than running a basic Bayesian update on the whole $T$ measurement set, which can be bounded as $O(T^3)$. 

\subsection{Role of the expert}

In SparQ-GP-UCB, the ``expert'' mechanism is not meant to be a human oracle, nor does it need to act as a perfectly accurate surrogate model. Instead, it serves as a means to partially refresh or correct stale information from previous observations in a principled and computationally bounded way.

More precisely, at each time step $t$, we are allowed to query the current value of the objective $f_t$ at a small number $\mathcal{O}\left( \log^d(t) \right)$ of previously observed points, selected via a $Q_t$-DPP sampling over $X_t$.

This mechanism is abstracted as an ``expert call'', but it is not assumed to be human or even a separate model. Rather, it reflects limited access to the current function values at previously observed locations, which can be interpreted in several realistic ways:

\begin{itemize}
    \item \textbf{Wireless sensor networks: } In Internet-of-Things applications~\cite{jamali2018internet}, sensors might collect data continuously but transmit selectively due to bandwidth or power constraints. Revisiting previous locations or reactivating a subset of sensors is often feasible, though costly — thus motivating a trade-off.
    \item \textbf{Physics-based monitoring: } In tasks such as environmental monitoring~\cite{roy2016spatio} where the underlying phenomenon is governed by a partial differential equation that needs to be simulated, the ``expert'' corresponds to access to the simulation itself. While running the simulator to evaluate the objective at a new point can be computationally expensive, it is often possible — though still costly — to re-run the simulator at previous input points to obtain updated objective values, reflecting changes in the underlying system.
    \item \textbf{Continual learning in ML systems: } For adaptive hyperparameter tuning or online systems~\cite{wang2024continual}, logs or cached evaluations might allow querying recent values again (e.g., checking performance of previous configurations on a new data batch).
\end{itemize}

We emphasize that the expert is not required to provide perfectly accurate information, but rather noisy or approximate values, consistent with a sub-Gaussian noise model. This is crucial in practice and aligns with many systems where re-evaluation is possible but noisy (e.g., due to changing conditions).

\section{Numerical results}

In this section, we compare the performance of SparQ-GP-UCB with four existing algorithms (TV-GP-UCB \cite{bogunovic2016}, W-GP-UCB \cite{deng2022weighted}, R-GP-UCB and SW-GP-UCB \cite{DBLP:Zhoujournals/corr/abs-2102-06296}) in a time-varying environment, on both a synthetic and a real-life dataset. We also run standard GP-UCB to show how it performs in time-varying settings. For all baseline methods, hyperparameters suchs as window size (R-GP-UCB, SW-GP-UCB), temporal kernel hyperparameter (TV-GP-UCB) and observations weights (W-GP-UCB) were set according to the recommendations provided by their respective authors.

\subsection{Synthetic data}

Observations are generated by perturbing the function 
\begin{align*}
  f\colon \Dset \times \mathbb{R}_+ & \to \mathbb{R}_+ \\[-1ex]
  (x,t) & \mapsto \exp(-0.05(x - 5\sin(0.1t))^2) + 0.5\cos(0.2 x) + 1.5
\end{align*}
with noise $ \epsilon \sim \mathcal{N}(0, \sigma^2)$, where the sampling noise variance $\sigma^2$ is set to $0.01$. We take the domain $\Dset=[-50,50]$. Let $t\geq0$. Then, by the mean value theorem 
$$
\underset{x \in \Dset}{\sup}\hspace{5pt}\vert f(x,t+1)-f(x,t)\vert \leq \underset{x,t}{\sup} \left\lvert \frac{\partial f(x,t)}{\partial t} \right\rvert \leq 1,
$$
and Assumption \ref{as.2} holds. We plot the average and standard deviation of the cumulative regret of each algorithm for $T=500$ iterations and $\delta = 0.05$ over 40 realizations using the squared exponential kernel, whose parameters have been fine-tuned by maximizing the log marginal likelihood of the data. For all four methods, we plot the mean and standard deviation of the average regret at each iteration. By selecting the number of queries $Q_T = 6\log(T)$ in line with the result of Proposition \ref{Thm:BoundKL}, we expect the average regret of SparQ-GP-UCB to vanish asymptotically. Furthermore, since the variation budget is not $V_T = o(T^{1/4})$ in our setting, due to the periodicity of $f \mapsto f(x,t)$, R-GP-UCB and SW-GP-UCB are not expected to have sublinear cumulative regret bounds. 

\begin{figure}[H]
\centering
    \includegraphics[scale=.45, trim=1cm 0cm 1cm 1cm, clip=on]{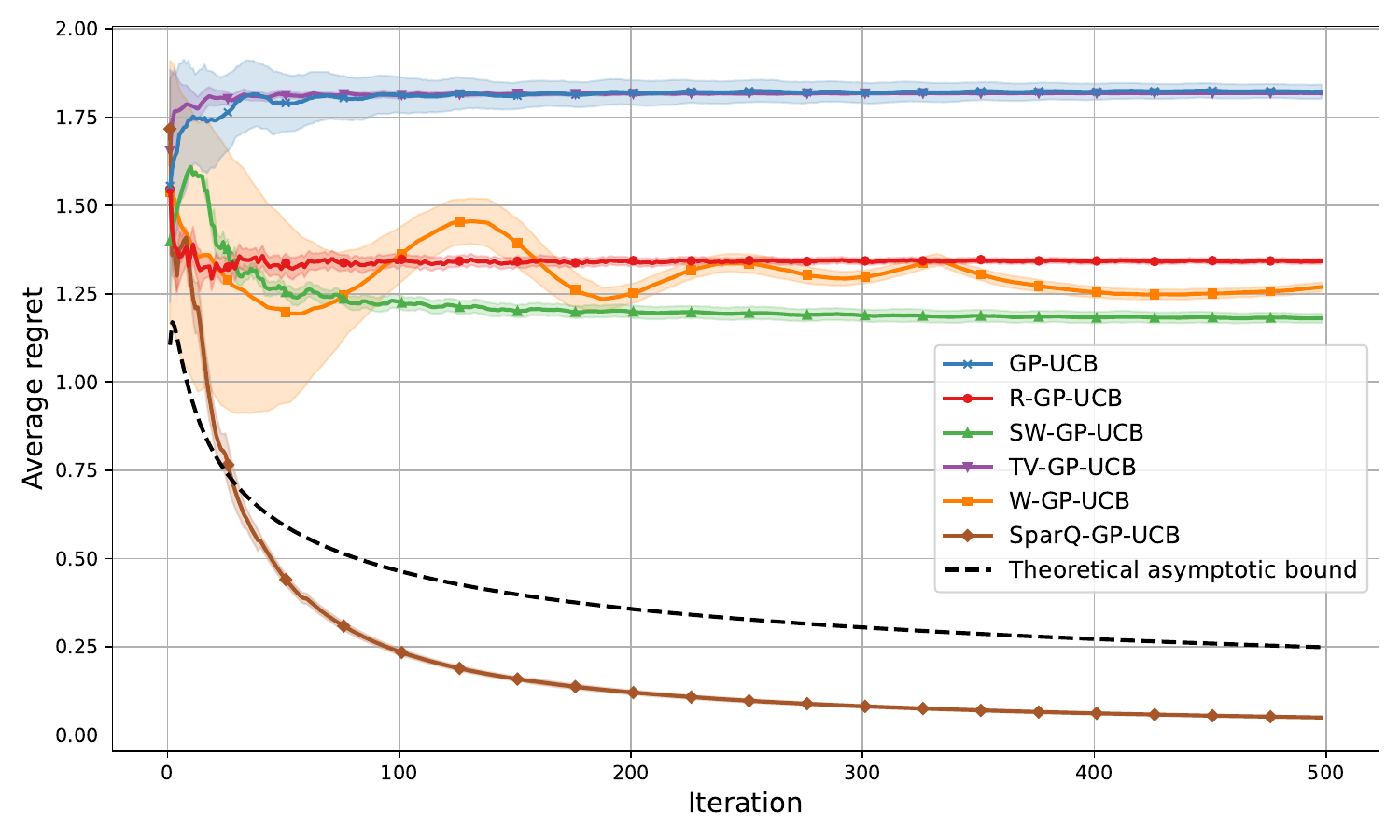}
    \caption{Average regret of GP-UCB variants in the time-varying setting.}
    \label{Fig:RegretStats}
\end{figure}

We can see in Figure~\ref{Fig:RegretStats} that SparQ-GP-UCB is the only method that converges to the dynamic optimum of the objective on average. Moreover, it falls below the theoretical bound (black curve) established in Eq. \eqref{Eq:RegretBound}. Standard TV methods (TV-GP-UCB, W-GP-UCB, R-GP-UCB and SW-GP-UCB) struggle to track the optimum and have linear cumulative regret ($R_T \approx 1.8$ for TV-GP-UCB, $R_T \approx 1.3$ for W-GP-UCB, $R_T \approx 1.34T$ for R-GP-UCB and $R_T \approx 1.18T$ for SW-GP-UCB). As expected, the average regret of standard GP-UCB grows slowly, suggesting a slightly superlinear average regret. In summary, at the cost of $\Ot(1)$ additional observations per iteration, SparQ-GP-UCB is the only method capable of accurately optimizing a time-varying objective with weak assumptions on its temporal variations.

\subsection{Real data}

To evaluate the effectiveness of SparQ-GP-UCB, we conducted experiments on a real-world dataset consisting of daily ozone level measurements collected from 28 sensors distributed across the New York City area over the course of two years \cite{epa2025dailydata}. This dataset presents a naturally time-varying environment, making it an ideal testbed for adaptive Bayesian optimization methods.

Again, we benchmarked SparQ-GP-UCB against several state-of-the-art baselines: GP-UCB, R-GP-UCB, SW-GP-UCB, TV-GP-UCB, and W-GP-UCB. The evaluation metric used was cumulative regret, plotted as average regret over time to highlight long-term performance trends. 

\begin{figure}[H]
\centering
    \includegraphics[scale=.45, trim=1cm 0cm 1cm 1cm, clip=on]{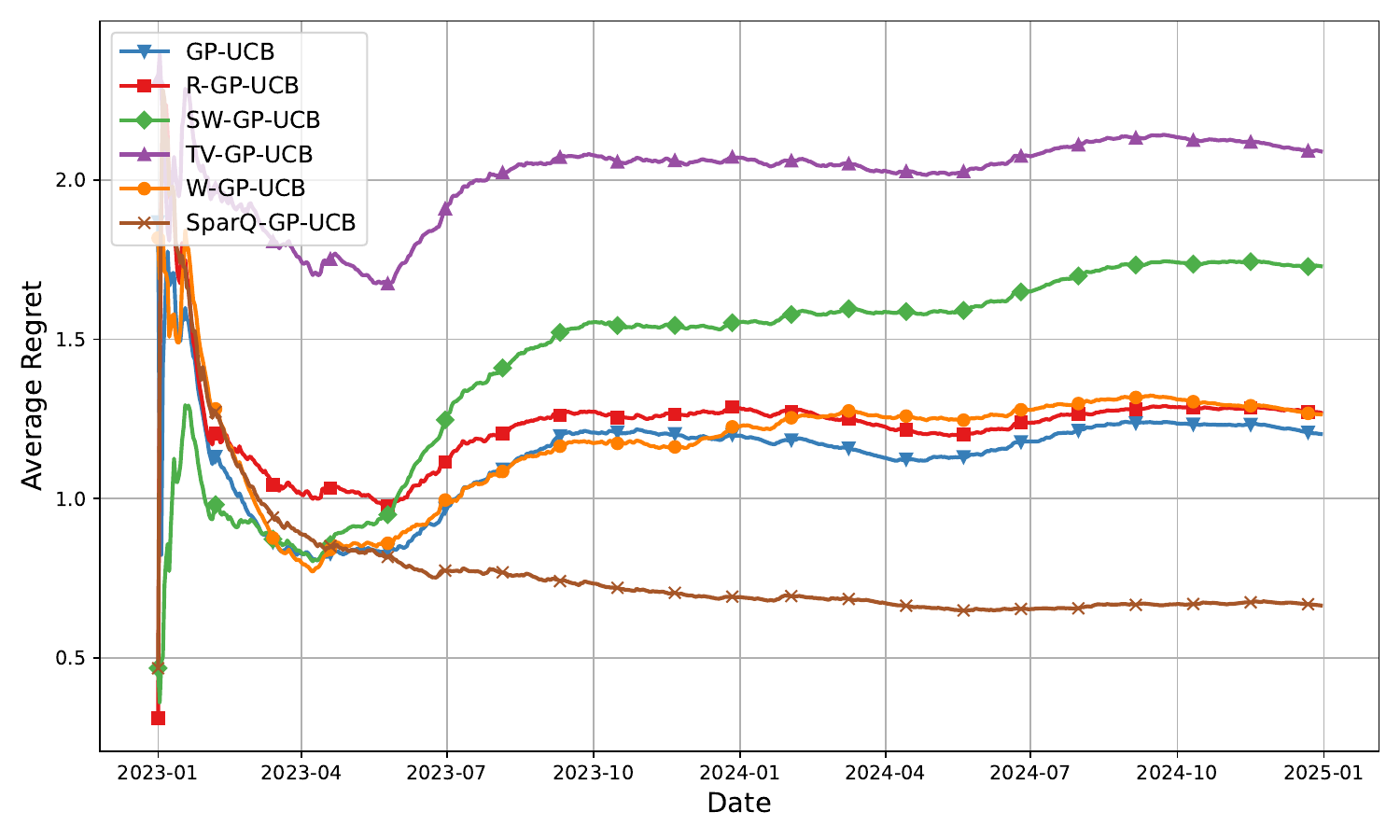}
    \caption{Numerical performances of GP-UCB variants on real data.}
    \label{Fig:RegretReal}
\end{figure}

Figure \ref{Fig:RegretReal} shows the evolution of average regret across time steps (here one day). We set $\delta=0.05$ and compute the kernel hyperparameters by maximizing the log marginal likelihood of the data. Although the ozone data is not generated by a model that respects our RKHS assumptions, we can see that SparQ-GP-UCB achieves significantly lower regret compared to other method. This indicates its superior ability to adapt to underlying non-stationary reward dynamics. Interestingly, GP-UCB achieves performance comparable to its time-varying counterparts. This may be attributed to characteristics of the ozone dataset — for example, the location of the maximum ozone level might be approximately stationary over time. Unlike the four time-varying baselines, SparQ-GP-UCB still consistently outperforms GP-UCB, suggesting that it effectively balances adaptation to temporal changes without overcompensating.

These results validate the robustness and adaptivity of our proposed method in capturing temporal variations and optimizing over dynamic environments.

\section{Proofs of main theorems and additional results}\label{sec:proofs}

\subsection{Proof of Theorem~\ref{Thm:regret}}

The regret proof is based on a few ingredients and extensions of previous work in~\cite{pmlr-v5-titsias09a,DBLP:journals/corr/abs-2111-03637,JMLR:v21:19-1015}. We proceed as follows: first we extend Proposition~\ref{Prop:PosteriorError} to the heteroscedastic and deterministic setting. Then we extend the regret results of~\cite{DBLP:journals/corr/abs-2111-03637} incorporating a sub-linearly growing maximum uncertainty, as well as the multiplicative error coming from the sparse inference. Then we combine the two results. 

It is convenient to define the set of pretend observations together with the latest ones with noise $< \sigma^2 \log(T)$ as $\Yv_T$ taken at points $\Xv_T$. The cardinality of the two sets is $T$. These sets are not the same as the set of sparse plus latest one observations, indicated as $\Ys_T$ and $\Xs_T$ whose cardinality is $\Ot(1)$. 

The first proposition extends the sparse approximation to our setting, and in particular, it is an extension of Corollary 22 in \cite{JMLR:v21:19-1015} to heteroscedastic GPs and deterministic inputs $\Xv_T$ in a compact domain.

\begin{proposition}\label{Thm:BoundKL}
    Let $\Xv_T$ be the set of actions chosen by SparQ-GP-UCB. Assume that pretend plus latest observations $\Yv_T \vert \Xv_T$ are conditionally Gaussian distributed. Then, under the same assumptions as Theorem \ref{Thm:regret}, for any $\eta >0$ and any $t$, there exists an approximation level $\varepsilon_t = \mathcal{O}\left( \frac{\eta}{t}\right)$ and number of queries $Q_t= \mathcal{O}\left( \log^d\left( \frac{t}{\eta} \right)\right)$ such that running SparQ-GP-UCB with an $\varepsilon_t$-approximate $Q_t$-DPP provides a posterior distribution $\Ps_T$ satisfying
    \begin{equation}\label{Eq:KLThm}
        \mathbb{E}[\mathrm{KL}[\Ps_T\|\Pv_T]] \leq \eta,
    \end{equation}
    where $\Pv_T$ is the posterior distribution on $(\Xv_T, \Yv_T)$, and $\mathrm{KL}$ is the KL divergence.
\end{proposition}


\begin{proof}
We give the proof in the Appendix. 
\end{proof}

While the result of Equation~\eqref{Eq:KLThm} is given in expectation, we can also use Markov's inequality implying,
\begin{equation}
\mathrm{KL}[\Ps_T\|\Pv_T] \leq \frac{2 \eta}{\delta},
\end{equation}
with probability $1-\delta/2$.

Let us now recall Proposition~1 of~\cite{JMLR:v21:19-1015}.

\begin{proposition}\label{Prop:PostError:app}[Proposition~1 of \cite{JMLR:v21:19-1015}]
    Let $P$ and $Q$ be the real and approximate posteriors with means $\mu_p$, $\mu_q$ and variances $\sigma_p^2$ and $\sigma_q^2$. Suppose $2\mathrm{KL}[Q\|P] \leq \eta \leq \frac{1}{5}$ and let $x \in \mathbb{R}^d$. Then,
    $$\vert \mu_p(x) - \mu_q(x) \vert \leq \sigma_p(x) \sqrt{\eta} \leq \frac{\sigma_q(x)\sqrt{\eta}}{\sqrt{1-\sqrt{3\eta}}} \hspace{10pt} \mbox{and} \hspace{10pt} \lvert 1 - \sigma_q^2/\sigma_p^2 \rvert < \sqrt{3\eta}.$$
\end{proposition}

Since by Proposition~\ref{Thm:BoundKL} we have a way to bound the error coming from considering a sparse setting instead of the pretend setting, and by Propostion~\ref{Prop:PostError:app}, we know how this translates into a multiplicative error of mean and variance, we are now ready for the regret result.  

\subsection{Regret proof}

\begin{proof}
In SparQ-GP-UCB algorithm, the posterior mean $\ms_t(.)$ and variance $(\sigs_t)^2(.)$ at step $t$ are obtained by performing regression on the sparse observations $(\Xs_{t-1}, \Ys_{t-1})$. The instantaneous regret of SparQ-GP-UBC at step $t$ is:
$$
r_t = f_t(x^*_t) - f_t(x_t),
$$
where, 
$$
x^*_t = \underset{x \in \Dset}{\text{argmax }} f_t(x) \hspace{10pt} \text{ and } \hspace{10pt} x_t = \underset{x \in \Dset}{\text{argmax }} \ms_{t-1}(x) + \beta_t \sigs_{t-1}(x).
$$
By leveraging the definition of confidence bounds acquisition functions $x \mapsto \text{ucb}_t(x)=\ms_{t-1}(x) + \beta_t \sigs_{t-1}(x)$ and $x \mapsto \text{lcb}_t(x) = \ms_{t-1}(x) - \beta_t \sigs_{t-1}(x)$, it is possible to bound the cumulative regret with probability $1-\delta/2$. To do that, we leverage the concentration bound provided in~\cite{kirschner2018informationdirectedsamplingbandits}.

\begin{proposition}\label{Prop:beta}{(Lemma~7, \cite{kirschner2018informationdirectedsamplingbandits})}
    Take any $0 < \delta \leq 1$ and let $f_T \in \mathcal{H}_k(\Dset)$ and $\mu_T(.)$ and $\sigma_T^2(.)$ be the posterior mean and covariance functions of $f_T(.)$ after observing $(X_T, Y_T)$ points. Then, for any $x \in \Dset$, the following holds with probability at least $1-\delta/2$:
    \begin{equation}\label{Eq:ucb}
        \forall t \in \{1, \hdots, T\}, \qquad \lvert \mu_{t-1}(x) - f_t(x) \rvert \leq \beta_t \sigma_{t-1}(x) 
    \end{equation}
where 
\begin{equation}\label{Eq:beta}
    \beta_t = \left( \sqrt{2 \log \left( \frac{2 \mathrm{det}\left( \Sigma_t + K_{tt} \right)^{1/2}}{\delta \mathrm{det}\left( \Sigma_t \right)^{1/2}} \right)} + \|f_t\|_k \right).
\end{equation}
\end{proposition}

Proposition is valid for the heteroscedastic setting. As such, with this in place, and with probability $1-\delta/2$: 
$$
r_t  \leq \text{ucb}_t(x^*_t) - \text{lcb}_t(x_t) \leq \text{ucb}_t(x_t) - \text{lcb}_t(x_t) = 2 \beta_t \sigs_{t-1}(x_t).
$$
Now we bound the cumulative regret at iteration $T$:
$$
R_T = \sum_{t=1}^T r_t \leq 2 \beta_T \sum_{t=1}^T \sigs_{t-1}(x_t).
$$
If we denote by $\sigv_{t-1}(.)$ the posterior variance of the regression on the pretend plus latest observations $(\Xv_{t-1}, \Yv_{t-1})$, Proposition \ref{Prop:PostError:app} gives
$$
\sigs_{t-1}(x_t) \leq \sigv_{t-1}(x_t) \sqrt{1+ \sqrt{3\eta}},
$$
with probability $1-\delta/2$, as long as we have a number of queries $\mathcal{O}(\log^d(t/\eta'))$ with $\eta' = \delta \eta/4$.

Therefore, the cumulative regret can be bounded with probability $1-\delta$ (for the union bound) as follows,
$$
R_T \leq 2 \beta_T \sqrt{1+ \sqrt{3\eta}} \sum_{t=1}^T \sigv_{t-1}(x_t).
$$
The observations $\Yv_{t-1}$ have been built such that their noise variance can be uniformly bounded by $\sigma^2 \log(t-1)$. By following the exact same computation steps of Makarova et al. in \cite{DBLP:journals/corr/abs-2111-03637} (Appendix A.1.1 Step 4) and replacing their fixed upper bound $\bar{\rho}$ by a logarithmically increasing upper bound $\sigma^2 \log(T)$ we get 
$$
R_T \leq 2 \beta_T \sqrt{1+ \sqrt{3\eta}} \sqrt{2T (1 + \left(\sigma^2\log(T)\right)^2) \gamma_T}\,,
$$
where $\gamma_T$ is the maximum information gain at step $T$. Finally,
\begin{equation}\label{Eq:OReg}
   R_T = \mathcal{O}\left(\beta_T \sqrt{T\log^2(T)\gamma_T}\right).
 \end{equation}

Let us now bound $\beta_T$ and $\gamma_T$. 

In SparQ-GP-UCB, the ucb acquisition function is computed using the approximate posterior mean and variance. We thus have: 
$$
\beta_T = \sqrt{2\log \left( \frac{2 \vert \textcolor{red}\Sigma^\mathrm{s}_{T} + \textcolor{red}{K}^\mathrm{s}_{TT} \vert^{1/2}} {\delta\vert \textcolor{red}\Sigma^\mathrm{s}_{T}\vert^{1/2}} \right)} + \|f_T\|_k.
$$

By the definition of information gain with the sparse plus recent observations (see, e.g.,~\cite{DBLP:journals/corr/abs-2111-03637}), we have
$$
\gamma_{Q_T} \geq \log \left(\frac{\vert \textcolor{red}\Sigma^\mathrm{s}_{T} + \textcolor{red}{K}^\mathrm{s}_{TT} \vert} {\vert \textcolor{red}\Sigma^\mathrm{s}_{T}\vert} \right),
$$
so that,
\begin{equation}\label{Eq:OBeta}
    \beta_T = \mathcal{O}\left(\sqrt{\log\left( \frac{2}{\delta} \right) + \gamma_{Q_T}}\right).
\end{equation}
If we combine bounds \eqref{Eq:OBeta} and \eqref{Eq:OReg}, we have a new expression for the regret bound:
\begin{equation}\label{Eq:OReg2}
    R_T = \mathcal{O}\left( \sqrt{\left(\log\left( \frac{2}{\delta} \right) + \gamma_{Q_T}\right)\left( T \log^2(T)\gamma_T\right)}\right).
\end{equation}

Again, by replacing $\bar{\rho}$ by $\sigma^ 2\log(T)$ in Makarova et al. proof (Appendix A.1.3) and using ${Q_T} = \mathcal{O}\left( \log^d\left( T\right)\right)$, we can bound the information gains $\gamma_T$ and $\gamma_{Q_T}$ for a squared exponential kernel\footnote{The information gain in a homoscedastic case for a SE kernel is $\mathcal{O}(d \log^{d+1}(T))$ to which we multiply a factor $\log^2(T)$ in our setting for the heteroscedastic case.}:
\begin{equation}\label{Eq:Gamma}
    \gamma_T = \mathcal{O} \left( d \log^{d+3}(T) \right),
\end{equation}
\begin{equation}\label{Eq:GammaStar}
        \gamma_{Q_T} = \mathcal{O} \left( d \log^{d+3}\left( \log^d \left(T \right) \right) \right)= \mathcal{O} \left( d \log^{d+3} \left(d\log(T) \right) \right).
\end{equation}

Finally, if we inject bounds \eqref{Eq:Gamma} and \eqref{Eq:GammaStar} into \eqref{Eq:OReg2}:
\begin{align}\label{Eq:OReg3}
    \begin{split}
        R_T &= \mathcal{O}\left( \sqrt{\left(\log\left( \frac{1}{\delta} \right) +  d \log^{d+3}\left(d \log \left(T\right) \right)\right)\left( T \log^2(T)d \log^{d+3}(T)\right)}\right)\\
        &= \mathcal{O}\left( \sqrt{\left(\log\left( \frac{1}{\delta} \right) +  d \log^{d+3}\left( d\log\left(T \right) \right)\right)\left( Td \log^{d+5}(T)\right)}\right) 
    \end{split}
\end{align}
This proves Theorem~\ref{Thm:regret}. \qed
\end{proof}

A closer look at the proof of~Theorem~\ref{Thm:regret} shows that one could choose to keep all the measurements with variance less than $g(T) = o(T^{1/4})$, as discussed in Section~\ref{sec:algo}, instead limiting at the ones with variance less than $\sigma^2\log(T)$. Since the maximum variance enters twice in the regret as a power of $2$, then the final regret would read $R = \Ot(\sqrt{T} \sqrt{g^4(T)}) = o(T)$, leading to a sublinear cumulative regret and a no-regret result.  

\subsection{Proof of Theorem~\ref{Thm:complexity}}

The computational  complexity of the algorithm proposed by Burt et al. \cite{JMLR:v21:19-1015} to obtain a $\varepsilon$ approximation of a $M$-DPP from a set of $N$ inputs is bounded as $\mathcal{O}\left(NM^3(\log\log N+\log M+\log 1/\varepsilon^2)\right)$, see their Section 4.2.2. 

The cost of the GP regression with $M$ training inputs is $\mathcal{O}(M^3)$ \cite{snelson2005sparse}, and the complexity of SparQ-GP-UCB is dominated by the computation of the $M$-DPP. Thus, for $T$ iterations in Algorithm~\ref{alg:SparQ-GP-UCB}  and with $Q_T$ the number of sparse inputs at the end of the process, the computational complexity of SparQ-GP-UCB is $T$ times the worst complexity of the DPP: 
$$
\mathcal{O}\left(T \left(Q_T\right)^3(\log \log T + \log Q_T + \log 1/\varepsilon_T^2)\right).
$$ 
In Proposition~\ref{Thm:BoundKL}, for fixed precision $\eta$, we show that $Q_T = \mathcal{O}\left(\log^d \left( T\right)\right)$, suffices to obtain a $\varepsilon_T$-approximation of a $Q_T$-DPP, with
$\varepsilon_T = \mathcal{O}\left( \frac{1}{T} \right)$. By substituting these estimates into the complexity, we obtain a total computational complexity of SparQ-GP-UCB of $\mathcal{O}\left(T^2 \log(T) \log^{3d}\left(T\right)\right)$. \qed

\section{Conclusion}

In this work, we provide a general framework to obtain sublinear regret bounds for GP optimization of a time-varying objective $f$ in the bandit setting. The function $f$ is assumed to belong to a RKHS with a bounded norm. We model time variations through uncertainty injection by linearly increasing the noise standard deviation of the data over time. We recover no-regret by asking $\Ot(1)$ additional side queries to an expert at each iteration. Future research will explore strategies to reduce the number of expert queries, such as retaining and reusing past responses to avoid querying the expert at every iteration.

\begin{credits}
\subsubsection{\ackname} This work was partly supported by the Agence Nationale de la Recherche (ANR) with the projects ANR AccelAILearning and ANR-23-CE48-0011-01.

\subsubsection{\discintname}
The authors have no competing interests to declare that are
relevant to the content of this article. 
\end{credits}

%
%
%

\newpage
\setcounter{page}{1}
\appendix

\begin{center}
\Large{\bf Supplementary material for: Time-varying Gaussian Process Bandit Optimization with Experts: no-regret in logarithmically-many side queries}
    
\end{center}

\section{Implementation details for the $M$-DPP algorithm }

In this section, we describe the algorithm mentioned in Proposition \ref{Prop:PosteriorError} and studied in the proof of Theorem \ref{Thm:complexity} used to select the $\Ot(1)$ sparse inputs, and especially a computationally lighter version introduced in \cite{JMLR:v21:19-1015} based on a Monte Carlo Markov Chain (MCMC) algorithm. A desirable property of sparse inputs is to be distributed in a way that captures information about the entire input set $X_T$. This can be done by leveraging the kernel matrix $K_T$ properties through the following process.

\begin{definition}{($M$-Determinantal Point Process)}
    Given a PSD kernel matrix $K$, a $M$-DPP is a discrete probability distribution defined over subsets of cardinality $M$ of the columns of $K$:
    $$\mathrm{Pr}(\mathbf{Z}=Z) = \frac{\mathrm{det}(K_{Z,Z})}{\sum_{\vert Z' \vert = M}\mathrm{det}(K_{Z',Z'})},$$
    where $K_{Z,Z}$ is the principal submatrix of $K$ with columns in $Z$.
\end{definition}
The probability of a set $Z$ is proportional to $\mathrm{det}(K_{Z,Z})$, which measures the diversity of the points in $Z$. As a consequence, $M$-DPP promotes dispersed subsets, which is in line with the intuition that the inducing inputs should be selected to guarantee good space coverage. 

However, obtaining an exact sample from an $M$-DPP, $M \in \mathbb{N}$, is polynomial in $M$ and nearly-linear in $N$ with $N$ the size of the dataset \cite{derezinski2019exact}. As the polynomial degree in $M$ is high, we use Algorithm 1 in Burt et al. \cite{JMLR:v21:19-1015} that draws approximate samples using a Markov chain algorithm. We report it in Algorithm~\ref{alg:M-DPP} with our notation. 

\begin{algorithm}[t]
\caption{MCMC algorithm to obtain an approximate sample from a $M$-DPP}\label{alg:M-DPP}
\begin{algorithmic}[1]
\Require Inputs $X=(x_i)_{i=1}^N$, number of sparse points $M$, kernel $k$, number of MCMC iterations $\tilde{N}$
\Ensure An approximate sample from an $M$-DPP: $Z_{\tilde{N}}$
\State Select an initial set of $M$ columns $Z_0$ to maximize the determinant of the associated submatrix
\For{$n \leq \tilde{N}$}
    \State Sample $i$ uniformly from $Z_n$ and $j$ uniformly from $X\setminus Z_n$
    \State Define $Z'=Z_n\setminus \{i\} \cup \{j\}$
    \State Compute $p_{i \rightarrow j}:= \frac{1}{2}\text{min}\left\{1, \frac{\text{det}\left(K_{Z'} \right)}{\text{det}\left(K_{Z_n} \right)}\right\}$
    \State With probability $p_{i \rightarrow j}$, $Z_{n+1}=Z'$ otherwise, $Z_{n+1}=Z_n$
\EndFor
\end{algorithmic}
\end{algorithm}

In order now to compute the first line of Algorithm~\ref{alg:M-DPP}, we require a bit more work. Let $C_M = (c_m)_{m=1}^M \subset X$, for $m \in \{1, \hdots, M\}$, $C_m=(c_i)_{i=1}^m \subset C_M$ and we denote by $\mathcal{S}_m \in \mathbb{R}$ the Schur complement of the principal submatrix $K_{C_{m-1}C_{m-1}}$ in 
$$K_{C_mC_m} = \begin{pmatrix} k(c_m,c_m) & k(c_m, C_{m-1}) \\ k(c_m, C_{m-1})^\top & K_{C_{m-1}C_{m-1}} \end{pmatrix}.$$
The expression of $\mathcal{S}_m$ is given by:
\begin{equation*}
    \mathcal{S}_m = k(c_m,c_m) - k(c_m, C_{m-1}) (K_{C_{m-1}C_{m-1}})^{-1} k(c_m, C_{m-1})^\top.
\end{equation*}
We have the following formula that enables to compute the determinant of $K_{C_mC_m}$ using $\text{det}(\mathcal{S}_m) = \mathcal{S}_m$ (because $\mathcal{S}_m \in \mathbb{R}$):
\begin{align}\label{Eq:SchurDeterminant}
\begin{split}
    \text{det}\left(K_{C_MC_M}\right) &= \mathcal{S}_M \text{det}\left(K_{C_{M-1}C_{M-1}}\right)\\
    &= \prod_{m=1}^M \mathcal{S}_m.
\end{split}
\end{align}
By leveraging Eq. \eqref{Eq:SchurDeterminant}, we build $Z_0$ by following the steps described in Algorithm \ref{alg:Z0}, where, for $x \in X$ and $Z \subset X$, $\mathcal{S}(x,Z)$ is the Schur complement of $K_{ZZ}$ in $K_{Z'Z'}$ where $Z' = Z \cup \{x\}$.

\begin{algorithm}[t]
\caption{Greedy sampling of $Z_0$}\label{alg:Z0}
\begin{algorithmic}[1]
\Require Inputs $X=(x_i)_{i=1}^N$, number of sparse points $M$, kernel $k$
\Ensure $Z_0$
\State $Z_0 = \emptyset$
\For{$m \leq M$}
    \State Sample $x_m = \underset{x \in X}{\text{argmax}}\hspace{3pt}\mathcal{S}(x,Z_0)$
    \State $Z_0 = Z_0 \cup \{x_m\}$
\EndFor
\end{algorithmic}
\end{algorithm}

\section{Proof of Proposition \ref{Thm:BoundKL} }

In order to prove Proposition \ref{Thm:BoundKL}, we start by showing two preliminary results. The first one provides an explicit expression for the KL divergence based on the expression of a lower bound for the true marginal likelihood. From this expression, we obtain an upper bound for $\mathrm{KL}[\Ps_T \| \Pv_T]$ that depends on the tail sum of the kernel matrix eigenvalues  $\Lambda_{Q_T}$ which leads us to demonstrate the second preliminary result yielding a bound for $\Lambda_{Q_T}$.

Note that for $\Xe \subset \Xv_T$ the set of sparse inputs selected and $\Ye \subset \Yv_T$ the set of observations delivered by the expert we have $\left(\Xe, \Ye\right) \subset (\Xs_T, \Ys_T)$ hence $\mathrm{KL}[\Ps_T \| \Pv_T] \leq \mathrm{KL}[\Pe \| \Pv_T]$, where $\Pe$ is the posterior distribution on $\left(\Xe, \Ye \right)$. In order to prove the result of Proposition \ref{Thm:BoundKL}, we settle for bounding $\mathrm{KL}[\Pe \| \Pv_T]$.

\subsection{Expression of the evidence lower bound}

Sparse inputs $\Xe$ and the distribution $\Phi$ of $f(\Xe) \vert \Yv_T$ are variational parameters selected to minimize the KL-divergence $\mathrm{KL}[\Pe \| \Pv_T]$. This minimization can equivalently be solved as the maximization of a lower bound for the true log marginal likelihood $\mathcal{L}(\Xe, \Phi)$ \cite{pmlr-v5-titsias09a} called the evidence lower bound (ELBO), since we have that (\cite{JMLR:v21:19-1015} Eq.~(8)),
\begin{equation}\label{Eq:ELBO}
    \mathcal{L}(\Xe, \Phi) + \mathrm{KL}[\Pe \| \Pv_T] = \log p(\Yv_T).
\end{equation}
In order to bound $\mathrm{KL}[\Pe \| \Pv_T]$, we start by explicitly expressing $\mathcal{L}(\Xe, \Phi)$.

For the sake of simplicity, we use generic notations in this section and show the following result, which is a generalization of a result from Titsias (Eq.~(9) in \cite{pmlr-v5-titsias09a}) to the heteroscedastic setting. 

\begin{proposition}\label{Prop:ELBO}
Consider $n$ observations $y$ at points $X$ with noise matrix $\Sigma$ and variational parameters $(X^*, p^*)$. After maximizing $\mathcal{L}$ with respect to the free variational distribution $p^*$, we get the following expression for the evidence lower bound: 
    \begin{equation}
        \mathcal{L}(X^*) = \log \mathcal{N}(y \vert 0, \Sigma + Q_{nn}) - \frac{1}{2}\mathrm{tr}(\Sigma^{-1}(K_{nn}-Q_{nn})),
    \end{equation}
    where $K_{nn}$ is the covariance matrix of the inputs $X$, $Q_{nn} = K_{n*}K_{**}^{-1}K_{n*}^\top$, $K_{n*}$ being the covariance matrix between $X$ and $X^*$ and $K_{**}$ the covariance matrix of $X^*$.
\end{proposition}

\begin{proof}
In order to prove these results, we follow the same steps as Titsias in the technical report \cite{Titsias2008VariationalMS} Section 3.1, with additional detailed explanations about some computations.
Let $\mathbf{f}$ be the true values of the reward at points $X$ and $P$ the exact posterior distribution. Sparse points are denoted as $(X^*, y^*)$, and $P^*$ is the approximate posterior distribution. We want to maximize the ELBO $\mathcal{L}(X^*, p^*)$ with respect to some free variational distribution $p^*$ to get a bound that only depends on $X^*$.

The joint density $p(y, \mathbf{f})$ is augmented with the variable $\mathbf{f}^*$ composed of the values of the reward evaluated at sparse inputs:
\begin{equation}\label{Eq:JointDensity}
    p(y, \mathbf{f}, \mathbf{f}^*) = p(y \vert \mathbf{f}) p(\mathbf{f} \vert \mathbf{f}^*)p(\mathbf{f}^*).
\end{equation}
By the marginalization property, we write the log marginal likelihood as
\begin{align*}
    \log p(y) &= \log \int \int p(y, \mathbf{f}, \mathbf{f}^*) \mathrm{d}\mathbf{f}\mathrm{d}\mathbf{f}^* \\
    &= \log \int \int p(y \vert \mathbf{f}) p(\mathbf{f} \vert \mathbf{f}^*)p(\mathbf{f}^*) \mathrm{d}\mathbf{f}\mathrm{d}\mathbf{f}^*.
\end{align*}
However, we also know that, 
\begin{align}\label{Eq:ELBOBis}
\begin{split}
    \log p(y) & = \mathcal{L}(X^*, p^*) + \mathrm{KL}[p^*(\mathbf{f}, \mathbf{f}^*) \| p(\mathbf{f}, \mathbf{f}^* \vert y)] \\
    & = \mathcal{L}(X^*, p^*) + \int_{\mathbf{f}^*} \int_{\mathbf{f}} p^*(\mathbf{f}, \mathbf{f}^*) \log \frac{p^*(\mathbf{f}, \mathbf{f}^*)}{p(\mathbf{f}, \mathbf{f}^* \vert y)} \mathrm{d}\mathbf{f} \mathrm{d}\mathbf{f}^*.
\end{split}
\end{align}
We use now Bayes' rule $p(\mathbf{f}, \mathbf{f}^* \vert y) = \frac{p(\mathbf{f}, \mathbf{f}^*, y)}{p(y)}$ to obtain a new expression for the KL divergence:

$$\mathrm{KL}(p^*(\mathbf{f}, \mathbf{f}^*) \| p(\mathbf{f}, \mathbf{f}^* \vert y)) = \underbrace{\int_{\mathbf{f}^*} \int_{\mathbf{f}} p^*(\mathbf{f}, \mathbf{f}^*) \log \frac{p^*(\mathbf{f}, \mathbf{f}^*)}{p(\mathbf{f}, \mathbf{f}^*, y)} \mathrm{d}\mathbf{f} \mathrm{d}\mathbf{f}^*}_{=-\mathcal{L}(X^*, p^*)} + \log p(y),$$
where we recognized one term as the ELBO. Finally, we inject Equation \eqref{Eq:JointDensity} in the expression of $\mathcal{L}(X^*, p^*)$ obtained above to get:

$$
\mathcal{L}(X^*, p^*) = \int_{\mathbf{f}^*} \int_{\mathbf{f}} p^*(\mathbf{f}, \mathbf{f}^*) \log \frac{p(y \vert \mathbf{f})p(\mathbf{f} \vert \mathbf{f}^*)p(\mathbf{f}^*)}{p^*(\mathbf{f}, \mathbf{f})} \mathrm{d}\mathbf{f} \mathrm{d}\mathbf{f}^*.
$$
As the exact posterior factorizes as $p(\mathbf{f}, \mathbf{f}^* \vert y) = p(\mathbf{f} \vert \mathbf{f}^*)p(\mathbf{f}^*,y)$, the variational distribution $p^*$ must satisfy the same factorization: $p^*(\mathbf{f}, \mathbf{f}^*) = p(\mathbf{f} \vert \mathbf{f}^*)\Phi(\mathbf{f}^*)$, where $\Phi$ is an \textit{unconstrained} variational distribution \cite{pmlr-v5-titsias09a}: 

\begin{align}\label{Eq:ELBOPhi}
\begin{split}
    \mathcal{L}(X^*, \Phi) &= \int_{\mathbf{f}^*} \int_{\mathbf{f}} p^*(\mathbf{f}, \mathbf{f}^*) \log \frac{p(y \vert \mathbf{f})p(\mathbf{f}^*)}{\Phi(\mathbf{f}^*)} \mathrm{d}\mathbf{f}\mathrm{d}\mathbf{f}^* \\
    & = \int_{\mathbf{f}^*} \Phi(\mathbf{f}^*) \Bigg\lbrace \underbrace{\int_{\mathbf{f}} \log (p(y\vert \mathbf{f})) p(\mathbf{f}\vert \mathbf{f}^*) \mathrm{d}\mathbf{f}}_{=\mathcal{I}(y,\mathbf{f})} + \log \frac{p(\mathbf{f}^*)}{\Phi (\mathbf{f}^*)} \Bigg\rbrace  \mathrm{d}\mathbf{f}^*.
\end{split}
\end{align}
Remember $y$ are noisy observations of the true values $\mathbf{f}$ of a function $f$ at points $X$, where the noise matrix is $\Sigma$. Then $y\vert \mathbf{f} \sim \mathcal{N}(\mathbf{f}, \Sigma)$ and 
%
%
\begin{equation}\label{Eq:PDFObs}
    \log(p(y \vert \mathbf{f})) = -\frac{n}{2} \log(2 \pi) - \frac{1}{2} \log (\vert \Sigma \vert ) - \frac{1}{2}(y-\mathbf{f})^\top \Sigma^{-1}(y-\mathbf{f}).
\end{equation}
Notice now that 
\begin{equation}\label{Eq:Rd1}
    (y-\mathbf{f})^\top \Sigma^{-1}(y-\mathbf{f}) = \mbox{tr}(\Sigma^{-1}(y-\mathbf{f})(y-\mathbf{f})^\top).
\end{equation}
By combining Eq.~\eqref{Eq:PDFObs} and Eq.~\eqref{Eq:Rd1}, we can get the following expression for $\mathcal{I}(y, \mathbf{f})$:
\begin{multline}\label{Eq:Int1}
    \mathcal{I}(y, \mathbf{f}) = -\frac{n}{2}\log (2 \pi) - \frac{1}{2} \log (\vert \Sigma \vert) \\ -\frac{1}{2} \int_{\mathbf{f}}p(\mathbf{f} \vert \mathbf{f}^*)\mbox{tr}(\Sigma^{-1}(yy^\top-2y\mathbf{f}^\top+\mathbf{f}\mathbf{f}^\top))\mathrm{d}\mathbf{f}.
\end{multline}   
If we now set $\alpha = \mathbb{E}[\mathbf{f}\vert\mathbf{f}^*]$ then 
\begin{multline*}
    \int_{\mathbf{f}}p(\mathbf{f} \vert \mathbf{f}^*)\mbox{tr}(\Sigma^{-1}(yy^\top-2y\mathbf{f}^\top+\mathbf{f}\mathbf{f}^\top)) \mathrm{d}\mathbf{f} = \int_{\mathbf{f}}p(\mathbf{f} \vert \mathbf{f}^*)\mbox{tr}\left(\Sigma^{-1}(yy^\top-2y\alpha^\top+\alpha\alpha^\top \right.\\\left.
    -2y(\mathbf{f}-\alpha)^\top  +2\alpha(\mathbf{f}-\alpha)^\top+
    (\mathbf{f}-\alpha)(\mathbf{f}-\alpha)^\top)\right)\mathrm{d}\mathbf{f}.
\end{multline*}
By the linearity of the trace operator, the right-hand side can be further bounded as,
\begin{multline*}
    \underbrace{\mbox{tr}(\Sigma^{-1}(yy^\top-2y\alpha^\top+\alpha\alpha^\top))}_{=(y-\alpha)^\top \Sigma^{-2}(y-\alpha)} - 2 \int_{\mathbf{f}}  p(\mathbf{f}\vert \mathbf{f}^*) \mbox{tr}(\Sigma^{-1}(y-\alpha)(\mathbf{f}-\alpha)^\top) \mathrm{d}\mathbf{f} + \\
    + \underbrace{\int_\mathbf{f} p(\mathbf{f}\vert \mathbf{f}^*) \mbox{tr}(\Sigma^{-1}(\mathbf{f}-\alpha)(\mathbf{f}-\alpha)^\top)\mathrm{d}\mathbf{f}}_{=\mbox{tr}(\Sigma^{-1}\mbox{Cov}(\mathbf{f}\vert \mathbf{f}^*))}.
\end{multline*}
But $\mathbf{f} \vert \mathbf{f}^* \sim \mathcal{N}(\alpha, K_{nn} - \underbrace{K_{n*}K_{**}^{-1}K_{n*}^\top}_{Q_{nn}})$ where $K_{nn}$ is the covariance matrix of the real inputs, $K_{**}$ the covariance matrix of the sparse inputs and $K_{n*}$ the covariance matrix between the real and sparse inputs so 
$$\mbox{tr}(\Sigma^{-1}\mbox{Cov}(\mathbf{f}\vert \mathbf{f}^*)) = \mathrm{tr}(\Sigma^{-1}(K_{nn}-Q_{nn})).$$
Moreover, 
\begin{align*}
    \int_{\mathbf{f}}  p(\mathbf{f}\vert \mathbf{f}^*) \mbox{tr}(\Sigma^{-1}(y-\alpha)(\mathbf{f}-\alpha)^\top) \mathrm{d}\mathbf{f} & = \mathbb{E}[\mathrm{tr}(\Sigma^{-1}(y-\alpha)(\mathbf{f}-\alpha)^\top) \vert \mathbf{f}^*] \\ &= 0.
\end{align*}
Therefore, by plugging in these expressions in Eq.~\eqref{Eq:Int1}, one gets:
$$
\mathcal{I}(y,\mathbf{f}) = \log [\mathcal{N}(y \vert \alpha , \Sigma^2)] - \frac{1}{2}\mathrm{tr}(\Sigma^{-1}(K_{nn}-Q_{nn})),
$$
and by using this result into Eq.~\eqref{Eq:ELBOPhi}, we obtain:
$$
\mathcal{L}(X^*, \Phi) = \int_{\mathbf{f}^*} \Phi(\mathbf{f}^*) \log \frac{\mathcal{N}(y \vert \alpha, \Sigma) p(\mathbf{f}^*)}{\Phi(\mathbf{f}^*)} \mathrm{d}\mathbf{f}^* - \frac{1}{2} \mathrm{tr}(\Sigma^{-1}(K_{nn}-Q_{nn})).
$$

We now want to maximize the evidence lower bound with respect to the distribution $\Phi$. As $\Phi$ is not constrained to belong to any restricted subset, the maximum is the equality case in Jensen's inequality: 
\begin{align*}
    \mathcal{L}(X^*) &= \log \left(\int_{\mathbf{f}^*} \mathcal{N}(y \vert \alpha, \Sigma)p(\mathbf{f}^*) \mathrm{d}\mathbf{f}^*\right) - \frac{1}{2} \mathrm{tr}(\Sigma^{-1}(K_{nn}-Q_{nn})) \\
    & = \log \left(\int_{\mathbf{f}^*} \mathcal{N}(y \vert \alpha, \Sigma)p(\mathbf{f}^*) \mathrm{d}\mathbf{f}^*\right) - \frac{1}{2} \mathrm{tr}(\Sigma^{-1}(K_{nn}-Q_{nn})) \\
    & = \log \mathcal{N}(y \vert 0, \Sigma + Q_{nn}) - \frac{1}{2}\mathrm{tr}(\Sigma^{-1}(K_{nn}-Q_{nn})).
\end{align*}
The last expression concludes the proof of Proposition~\ref{Prop:ELBO}. \qed
\end{proof}

Given Equation \eqref{Eq:ELBOBis} and Proposition \ref{Prop:ELBO} we can write the KL divergence between the sparse and pretend posteriors as:

\begin{align}\label{Eq:KLTrace}
\begin{split}
    \mathrm{KL}[\Pe\|\Pv_T] &= \log p(\Yv_T) -  \mathcal{L}\left(\Xe\right) \\
    &= \log \frac{\mathcal{N}(\Yv_T \vert 0,\Sigv_{TT} + K_{TT})}{\mathcal{N}(\Yv_T \vert 0, \Sigv_{TT} + Q_{TT})} + \frac{1}{2}\mbox{Tr}((\Sigv_{TT})^{-1}(K_{TT}-Q_{TT})),
\end{split}
\end{align}
where $Q_{TT} = K_{T\mathrm{v}}(K^\mathrm{E}_{TT})^{-1} K_{T\mathrm{v}}^\top$, $K_{T\mathrm{v}}$ being the covariance matrix between $\Xv_T$ and $\Xe$ and $K^\mathrm{E}_{TT}$ is the covariance matrix of $\Xe$.

\subsection{Bounding the tail sum of the kernel eigenvalues}

Let $\tilde{\lambda}_1 \geq \hdots \geq \tilde{\lambda}_T \geq 0$ be the eigenvalues of the kernel matrix $K_{TT}$. The objective of this section is to bound 
$$
\Lambda_{Q_T} = \sum_{t=Q_T+1}^T \tilde{\lambda}_t.
$$
To do that, we use the spectrum of the kernel operator associated with some density $q(.)$ according to measure $\mu$, $\mathcal{K}_q : L^2(\mathbb{R}^d, \mu) \to L^2(\mathbb{R}^d, \mu)$ defined as follows (Koltchinskii and Giné, 2000) \cite{bj/1082665383}:
\begin{equation}\label{Eq:KernelOpeator}
   (\mathcal{K}_qg)(x) = \int g(x')k(x,x')p(x')\mathrm{d}x'. 
\end{equation}
Bounds for $\Lambda_{Q_T}$ have been obtained for stochastic covariates $X_T$ and, to our knowledge, we are the first to extend this result to deterministic datasets.

\begin{theorem}\label{Thm:TailSum}
    Let $\Dset \subset \mathbb{R}^d$ be convex and compact. Suppose that $\Xv_T = (x_t)_{t=1}^T \subset \Dset$ is the set of inputs obtained after $T$ iterations of SparQ-GP-UCB. Let $q(.)$ be the PDF of the Gaussian distribution $\mathcal{N}(\mu_q,l_q/8)$, where $\mu_q = \arg \underset{x \in \mathbb{R}^d}{\min} \underset{y \in \Dset}{\max} \|x-y\|$ is the Chebyshev center of $\Dset$ and $l_q = \underset{x,y\in\Dset}{\text{max}} \hspace{5pt}\|x-y\|$. Let $\lambda_1 \geq \lambda_2 \geq\hdots $ be the eigenvalues of the integral operator $\mathcal{K}_q$. Then, for all $Q_T\geq 1$, we have that
    \begin{equation}
        \mathbb{E}\left[\frac{\Lambda_{Q_T}}{T} \right] \leq C_{\Dset} \sum_{t=Q_T+1}^{\infty} \lambda_t,
    \end{equation}
    where $C_{\Dset} = \left( \sqrt{2\pi} \frac{l_q}{8}\right)^{d} e^8$.
\end{theorem}

Before studying the kernel behavior for deterministic inputs and therefore providing the proof, we provide some background about how to bound $\Lambda_{Q_T}$ when inputs are stochastic.

\smallskip
{\bf Stochastic inputs.} For stochastic i.i.d. covariates with continuous density $p$ according to some measure $\mu$, the matrix $K_{TT}/T$ behaves like the operator $\mathcal{K}_p$. Following this result, Shawe-Taylor et al. \cite{shawe2005} provide a bound for the tail sum $\Lambda_{Q_T}$ and Burt et al. \cite{JMLR:v21:19-1015} generalize their result even when the covariates are not identically distributed (Lemma 11). We slightly adapt their notation for compactly supported covariates in the following lemma.

\begin{lemma}\label{Lemma:TailBound}
    Suppose that each of the the $T$ covariates $(\mathbf{x}_t)_{t=1}^T$ has compact support $\Dset$, continuous density $p_t(.) : \Dset \to [0,1]$ and that there exists a density $q(.)$ and constants $(c_t)_{t=1}^T$ such that for all $x \in \mathbb{R}^d$ and for all $t \in \{1, \hdots, T\}$, $p_t(x) \leq c_t q(x)$. Let $\lambda_1 \geq \lambda_2 \geq\hdots $ be the eigenvalues of the integral operator $\mathcal{K}_q$. Then, for all $Q_T\geq 1$
    $$\mathbb{E}\left[\frac{\Lambda_{Q_T}}{T} \right] \leq \bar{c} \sum_{t=Q_T+1}^{\infty} \lambda_t$$
    where $\bar{c} = \frac{1}{T} \sum_{t=1}^T c_t$. 
\end{lemma}

It is possible to evaluate the sum $\sum_{t={Q_T}+1}^{\infty} \lambda_t$ for specific kernels and distributions. In particular, we have the following result for $q(.)$ Gaussian and $k(.,.)$ the squared exponential (SE) kernel.

\begin{proposition}\label{Prop:KernelOpBound}[Proposition 21 of \cite{JMLR:v21:19-1015}]
    If $k(.,.)$ is the SE-kernel and $q(.)$ a Gaussian density, then, for $Q_T \geq \frac{1}{\alpha}d^d + d - 1$, there exists $C > 0$ such that $\sum_{t=Q_T+1}^{\infty} \lambda_t \leq CQ_T \exp(-\alpha (Q_T)^{1/d})$ where $\alpha > 0$ depends on the model parameters.
\end{proposition}

It is important to note that, in our configuration, Lemma \ref{Lemma:TailBound} and Proposition \ref{Prop:KernelOpBound} cannot be applied ``as is'' due to the fact that the dataset $\Xv_T = (x_t)_{t=1}^T$ is sequentially deterministic. In fact, each $x_i$ is deterministic knowing $(x_1, \hdots, x_{i-1})$, and $x_1$ is deterministic given the prior $\mathcal{N}(\mu_0(x) = 0, \sigma^2_0(x)=k(x,x))$. Indeed, for $t \in \{1, \hdots, T\}$, $x_t$ is obtained by solving the optimization problem,
$$
x_t = \underset{x \in \Dset}{\text{argmax}} \hspace{5pt} \underbrace{\ms_{t-1}(x) + \beta_t \sigs_{t-1}(x)}_{=\text{ucb}_t(x)}. 
$$

\smallskip
{\bf Deterministic inputs and proof of Theorem \ref{Thm:TailSum}.} To tackle deterministic inputs, the idea is to see $\Xv_T$ as a non-Markovian process where the probability of each event depends on the state of all previous events. At time $t \in \{1, \hdots, T\}$, we view the UCB acquisition function $x \mapsto \text{ucb}_t(x)$ as defining a likelihood over potential values of $\Dset$ and define the probability density function: 
$$
p_t(x) = \mathds{1}_{\Dset}(x) \frac{\exp\left(\frac{\text{ucb}_t(x')}{\tau}\right)}{\int_{\Dset}\exp \left(\frac{\text{ucb}_t(x')}{\tau}\right)\mathrm{d}x'},
$$
where $\mathds{1}_{\Dset}$ is the indicator function of $\Dset$. The density function $p_t(.)$ is compactly supported (its support is the compact set $\Dset$) and $\int_{\mathbb{R}^d}p_t(x)\mathrm{d}x = 1$. Parameter $\tau > 0$ controls the sharpness of the distribution and draws the link between the deterministic and probabilistic settings as 
$$
p_t(x) \underset{\tau \to \infty}{\longrightarrow} \delta (x-x_t),
$$
where $\delta (.)$ is the Dirac distribution.

Let $q(.)$ be the density of the Gaussian distribution $\mathcal{N}(\mu_q, l_q/8)$ where $\mu_q$ is the center of $\Dset$ and $l_q = \underset{x , y \in \Dset}{\text{max}} \hspace{5pt} \|x-y\|$. Then, for $x \in \Dset$, 
\begin{equation*}
    \frac{p_t(x)}{q(x)} \leq \frac{1}{q(x)} = \left(2 \pi \left(\frac{l_q}{8}\right)^2\right)^{d/2} \exp\left( \frac{64\|x-\mu_q\|^2}{2l_q^2} \right) \leq \left( \sqrt{2\pi} \frac{l_q}{8}\right)^{d} e^8.
\end{equation*}
If we set $c=\left( \sqrt{2\pi} \frac{l_q}{8}\right)^{d} e^8$, then for all $t=1, \hdots, T$ and for all $x \in \mathbb{R}^d$, $p_t(x) \leq cq(x)$. As a conclusion, each $x_t \in \Xv_T$ can be seen as an observation of a covariate $\mathbf{x}_t$ with continuous density $p_t(.) \leq cq(.)$, where $q(.)$ is the density of a multivariate Gaussian distribution. Lemma \ref{Lemma:TailBound} can thus be applied and we obtain:
\begin{equation}\label{Eq:TailBound}
    \mathbb{E}[\Lambda_{Q_T}] \leq C_{\Dset}T \sum_{t=Q_T+1}^{\infty} \lambda_t.
\end{equation}
This proves Theorem \ref{Thm:TailSum}. \qed

Now if we take $\alpha$, $C$ and $Q_T$ as in Proposition \ref{Prop:KernelOpBound} and apply the result to inequality \eqref{Eq:TailBound}, we obtain
\begin{equation}\label{Eq:TailBound2}
    \mathbb{E}[\Lambda_{Q_T}] \leq \bar{C} Q_TT \exp\left(-\alpha (Q_T)^{1/d}\right),
\end{equation}
where $\bar{C} = C_{\Dset}C$.

\subsection{Bounding the KL divergence in probability and proof of Proposition \ref{Thm:BoundKL}}

Now that we have shown these two preliminary results, we have the material necessary to complete the proof of Proposition \ref{Thm:BoundKL}. Equation \eqref{Eq:KLTrace} gives
\begin{equation}\label{Eq:KLTrace:bis}
    \mathrm{KL}[\Pe\|\Pv_T] = \log \frac{\mathcal{N}(\Yv_T \vert 0,\Sigv_{TT} + K_{TT})}{\mathcal{N}(\Yv_T \vert 0, \Sigv_{TT} + Q_{TT})} + \frac{1}{2}\mbox{Tr}((\Sigv_{TT})^{-1}(K_{TT}-Q_{TT})).
\end{equation}
We now take the conditional expectation $\mathbb{E}[. \vert {\textcolor{blue}{\mathbf{X}}^\mathrm{v}} = \Xv_T, {\mathbf{E}}^\mathrm{E} = \Xe]$:
\begin{multline}\label{Eq:Ep}
    \mathbb{E}[\mathrm{KL}[\Pe\|\Pv_T] \vert{\textcolor{blue}{\mathbf{X}}^\mathrm{v}} = \Xv_T, {\mathbf{E}}^\mathrm{E} = \Xe] =\\
    \mathbb{E}\left[\log \frac{\mathcal{N}(\Yv_T \vert 0,\Sigv_{TT} + K_{TT})}{\mathcal{N}(\Yv_T \vert 0, \Sigv_{TT} + Q_{TT})} \vert {\textcolor{blue}{\mathbf{X}}^\mathrm{v}} = \Xv_T, {\mathbf{E}}^\mathrm{E} = \Xe\right] + \frac{1}{2}\mbox{Tr}((\Sigv_{TT})^{-1}(K_{TT}-Q_{TT})) \\
     = \underbrace{\mathrm{KL}[\mathcal{N}(\Yv_T \vert 0, \Sigv_{TT}+ K_{TT}) \| \mathcal{N}(\Yv_T \vert 0, \Sigv_{TT} + Q_{TT})]}_{=\textsf{kl}(K_{TT},Q_{TT})} +\\+ \frac{1}{2}\mbox{Tr}((\Sigv_{TT})^{-1}(K_{TT}-Q_{TT}))
    \leq \textsf{kl}(K_{TT},Q_{TT}) + \frac{1}{2\lambda_{\min}(\Sigv_{TT})}\mbox{Tr}(K_{TT}-Q_{TT}).
\end{multline}

Following the proof of Lemma~4 of \cite{JMLR:v21:19-1015}, we bound the term $\textsf{kl}(K_{TT},Q_{TT})$:
\begin{align*}
    \textsf{kl}(K_{TT},Q_{TT}) &= \frac{1}{2}\Big( \log(\vert \Sigv_{TT} + Q_{TT}\vert ) - \log(\vert \Sigv_{TT} + K_{TT}\vert ) - T + \\ & \hspace*{2cm} \text{tr}((\Sigv_{TT} + Q_{TT})^{-1}(\Sigv_{TT} + K_{TT})) \Big) \\
    & \leq \frac{1}{2}\left( - T + \text{tr}((\Sigv_{TT} + Q_{TT})^{-1}(\Sigv_{TT} + K_{TT})) \right)\\
    &= \frac{1}{2} \left( - T + \text{tr}((\Sigv_{TT} + Q_{TT})^{-1}((\Sigv_{TT} + Q_{TT})+(K_{TT}-Q_{TT})) \right) \\
    &= \frac{1}{2}\text{tr}\left( (Q_{TT}+\Sigv_{TT})^{-1}(K_{TT}-Q_{TT}) \right),
\end{align*}
where the inequality in the second line is due to the fact that $Q_{TT}$ is a rank $Q_T < T$ approximation of $K_{TT}$. We use a special case of Hölder inequality (Proposition~36 in Burt et al. \cite{JMLR:v21:19-1015}):
\begin{proposition}[Proposition~36 of \cite{JMLR:v21:19-1015}]\label{Prop:Holder}
    Let $A$ and $B$ be symmetric positive semi-definite matrices. Then:
    \begin{equation}\label{Eq:Holder}
        \mathrm{tr}\left(BA\right) \leq \mathrm{tr}(A)\lambda_{\max}(B).
    \end{equation}
\end{proposition}
We apply Equation \eqref{Eq:Holder} with $A = K_{TT}-Q_{TT}$ and $B = (Q_{TT}+\Sigv_{TT})^{-1}$: 
\begin{equation}\label{Eq:kl}
\textsf{kl}(K_{TT},Q_{TT}) \leq \frac{1}{2\lambda_{\min}(\Sigv_{TT})}\mbox{tr}(K_{TT}-Q_{TT}).
\end{equation}
Finally, we combine Equations \eqref{Eq:Ep} and \eqref{Eq:kl} to obtain:
\begin{equation*}
    \mathbb{E}[\mathrm{KL}[\Pe_T\|\Pv_T] \vert {\textcolor{blue}{\mathbf{X}}^\mathrm{v}} = \Xv_T, {\mathbf{E}}^\mathrm{E} = \Xe] \leq \frac{1}{\lambda_{\min}(\Sigv_{TT})}\mbox{tr}(K_{TT}-Q_{TT}).
\end{equation*}
We recall that $\Sigv_{TT}$ is the noise matrix of pretend plus latest observations hence $\lambda_{\min}(\Sigv_{TT}) = \sigma^2$ and finally:

\begin{equation}\label{Eq:BoundKL}
    \mathbb{E}[\mathrm{KL}[\Pe_T\|\Pv_T] \vert {\textcolor{blue}{\mathbf{X}}^\mathrm{v}} = \Xv_T, {\mathbf{E}}^\mathrm{E} = \Xe] \leq \frac{\mbox{tr}(K_{TT}-Q_{TT})}{\sigma^2}.
\end{equation}

Let $\varepsilon_T > 0$ be some error. Sparse inputs $\Xs_T$ are subsampled from $\Xv_T$ according to a $\varepsilon_T$-approximation of a $Q_T$-DPP (Algorithm 1 in \cite{JMLR:v21:19-1015}) and we use the following proposition to bound the trace (Corollary 8 in \cite{JMLR:v21:19-1015}).

\begin{proposition}\label{Cor:MDPP}[Corollary~8, \cite{JMLR:v21:19-1015}]
    Let $\rho$ denote a $M$-DPP with kernel matrix $L \in \mathbb{R}^{N \times N}$, satisfying $L_{i,i} \leq v$. Let $\rho'$ denote a measure on subsets of columns of $L$ with cardinality $M$ such that $\mathrm{TV}(\rho, \rho') \leq \varepsilon_T$, where the total variation is defined as: $\mathrm{TV}(\rho, \rho'):= \frac{1}{2}\sum_{\vert Z \vert = M} \rho(Z) - \rho'(Z)$. Then 
    $$
    \mathbb{E}_{\rho'}[\mathrm{Tr}(L - L_Z)] \leq 2Nv\varepsilon_T + (M+1)\sum_{m=M+1}^N \eta_m,
    $$
    where $\eta_m$ is the $m^{th}$ largest eigenvalue of $L$ and $L_Z = L_{Z,N}^T L_{Z,Z}^{-1}L_{Z,N}$.
\end{proposition}
Applying Proposition \ref{Cor:MDPP} to our setting, by noticing that Algorithm 1 in \cite{JMLR:v21:19-1015} provides a way to have $\mathrm{TV}(\rho, \rho') \leq \varepsilon_T$ (see Lemma~9 of \cite{JMLR:v21:19-1015}), and by setting $L = K_{TT}$, $L_Z = Q_{TT}$ and $v = M_k^2$, we have
\begin{equation}\label{Eq:Trace}
    \mathbb{E}_{\rho '}[\mathrm{Tr}(K_{TT}-Q_{TT})] \leq 2TM_k^2 \varepsilon_T + (Q_T+1)\sum_{t=Q_T+1}^T \tilde{\lambda}_t,
\end{equation}
where $\tilde{\lambda}_1 \geq \hdots \geq \tilde{\lambda}_T$ are the eigenvalues of the kernel matrix $K_{TT}$. We can see that the quality of the approximation depends on the decay rate of the kernel spectrum. The faster the decay, the fewer points are needed to accurately approximate the virtual distribution. 

Taking the expectation of Equation \eqref{Eq:BoundKL} with respect to the sparse input distribution over the sparse inputs distribution we get:
\begin{equation}\label{Eq:BoundKLExp}
    \mathbb{E}[\mathrm{KL}[\Pe_T\|\Pv_T] \vert {\textcolor{blue}{\mathbf{X}}^\mathrm{v}} = \Xv_T] \leq \frac{\mathbb{E}_{\rho'}[\mbox{tr}(K_{TT}-Q_{TT})]}{\sigma^2}.
\end{equation}
From Equations \eqref{Eq:BoundKLExp} and \eqref{Eq:Trace}, we get a bound for the total expectation of the KL-divergence conditioned on $\textcolor{blue}{\mathbf{X}}^\mathrm{v} = \Xv_T$ as,
\begin{equation}\label{Eq:BoundKLExp2}
 \mathbb{E}[\mathrm{KL}[\Pe_T\|\Pv_T]\vert {\textcolor{blue}{\mathbf{X}}^\mathrm{v}} = \Xv_T] \leq  \frac{2TM_k^2 \varepsilon_T + (Q_T+1)\Lambda_{Q_T}}{\sigma^2}.
\end{equation}
Now, taking the expectation with respect to the covariate distribution over the covariate distribution of Equation \eqref{Eq:BoundKLExp2} and applying Equation \eqref{Eq:TailBound2} we obtain:
\begin{equation*}
 \mathbb{E}[\mathrm{KL}[\Pe_T\|\Pv_T]] \leq  \frac{2TM_k^2 \varepsilon_T + (Q_T+1)\bar{C} Q_TT \exp\left(-\alpha (Q_T)^{1/d}\right)}{\sigma^2}.
\end{equation*}

In order to find appropriate $\varepsilon_T$ and $Q_T$ such that 
$$
\mathbb{E}[\mathrm{KL}[\Pe_T\|\Pv_T]] \leq \eta,
$$
we solve the system of equations
\begin{subequations}
\begin{align}
\frac{2TM_k^2\varepsilon_T}{\sigma^2} &= \frac{\eta}{2}, \label{eqn:line-1} \\
\frac{\bar{C}TQ_T(Q_T+1)\exp\left(-\alpha (Q_T)^{1/d}\right)}{\sigma^2}
 &\leq \frac{\eta}{2}. \label{eqn:line-2}
\end{align}
\label{eqn:all-lines}
\end{subequations}
Equation~\eqref{eqn:line-1} gives
\begin{equation}\label{Eq:epsilon}
    \varepsilon_T = \frac{\sigma^2 }{4TM_k^2}\eta = \mathcal{O}\left( \frac{\eta}{T}\right).
\end{equation}

It is reasonable to assume $Q_T<T$, otherwise $\Xs_T=X_T$ and the sparse and pretend posteriors coincide, which would imply
$$\mathbb{E}[\mathrm{KL}[\Pe_T\|\Pv_T]] = 0.$$
For $Q_T<T$, a valid solution for Equation~\eqref{eqn:line-2} is given by the solution of:
\begin{align}\label{Eq:SparseTime}
\begin{split}
  & \frac{\bar{C}TQ_T(Q_T+1)\exp\left(-\alpha (Q_T)^{1/d}\right)}{\sigma^2} < \frac{\bar{C}T^3\exp\left(-\alpha (Q_T)^{1/d}\right)}{\sigma^2} = \frac{\eta}{2}\\
  \Leftrightarrow \qquad & \exp\left(\alpha Q_T^{1/d}\right) = \frac{2\bar{C}T^3}{\sigma^2 \eta}\\
  \Leftrightarrow \qquad &Q_T= \left( \frac{1}{\alpha} \right)^d \log^d\left( \frac{2 \bar{C}T^3}{\sigma^2 \eta} \right) =\mathcal{O}\left( \log^d(T/\eta)\right).
\end{split}
\end{align}
By taking $\epsilon_T$ and $Q_T$ according to \eqref{Eq:epsilon} and \eqref{Eq:SparseTime}, we get 
$$\mathbb{E}[\mathrm{KL}[\Pe_T\|\Pv_T]] \leq \eta,$$
and we have proven Proposition \ref{Thm:BoundKL}. \qed

This proposition states that we can get an approximate posterior $\Ps_T$ arbitrarily close to the virtual posterior $\Pv_T$ in terms of KL-divergence by running SparQ-GP-UCB on our dataset with appropriate choices for $\varepsilon_T$ and $Q_T$.

\end{document}